\documentclass[11pt]{amsart}
\usepackage{amssymb,latexsym,amsmath,graphicx,graphics,epic,eepic}

\theoremstyle{plain}
\newtheorem{theorem}{Theorem}
\numberwithin{theorem}{section}
\newtheorem{lemma}[theorem]{Lemma}
\newtheorem{proposition}[theorem]{Proposition}
\newtheorem{corollary}[theorem]{Corollary}

\theoremstyle{definition}

\newtheorem{example}[theorem]{Example}

\newtheorem{remark}[theorem]{Remark}

\newcommand{\C}{{\mathbb C}}
\newcommand{\R}{{\mathbb R}}
\newcommand{\Z}{{\mathbb Z}}
\newcommand{\Q}{{\mathbb Q}}
\renewcommand{\P}{{\mathbb P}}
\newcommand{\s}{{\mathbb S}}

\newcommand{\I}{{\mathbb I}}

              \newcommand{\M}{{\mathcal M}}

\begin{document}
\title{Symplectic symmetries of $4$-manifolds} 
\author{Weimin Chen and Slawomir Kwasik}
\thanks{The first author was supported in part by NSF grant
DMS-0304956.}
\date{}

\begin{abstract}
A study of symplectic actions of a finite group $G$ on smooth $4$-manifolds
is initiated. The central new idea is the use of $G$-equivariant
Seiberg-Witten-Taubes theory in studying the structure of the fixed-point
set of these symmetries. The main result in this paper is a complete 
description of the fixed-point set structure (and the action around it) of
a symplectic cyclic action of prime order on a minimal symplectic 
$4$-manifold with $c_1^2=0$. Comparison of this result with the case of 
locally linear topological actions is made. As an application of these 
considerations, the triviality of many such actions on a large class
of $4$-manifolds is established. In particular, we show the triviality
of homologically trivial symplectic symmetries of a $K3$ surface (in
analogy with holomorphic automorphisms). Various examples and comments 
illustrating our considerations are also included.
\end{abstract}

\maketitle

\section{Introduction}

When studying smooth or locally linear, topological actions of a finite
group $G$ on $4$-manifolds the central problem is to describe the 
structure of the fixed-point set and the action around it. This together
with the $G$-signature theorem of Atiyah-Singer (cf. \cite{AS, HZ})
leads to a wealth of information about the action of $G$. For smooth
actions, there is additional information provided through the use of gauge
theory (cf. e.g. \cite{D1, DK, Fu, BF}). Indeed, the $G$-signature 
theorem (and 
gauge theory in the case of smooth actions) imposes restrictions on the
symmetries of $4$-manifolds through the fixed-point set and the action
around it. However, in general no substantial insight can be obtained without
more specific knowledge of the latter. In fact, Freedman's $4$-dimensional
topological surgery theory (cf. \cite{F,FQ}) allows various constructions
of periodic homeomorphisms on simply-connected $4$-manifolds. In particular, 
in the case of topological, locally linear actions results of Edmonds
(cf. \cite{E}) demonstrate great flexibility of such actions, and the work of
Edmonds and Ewing (cf. \cite{EE}) shows that for locally linear, pseudofree
(i.e. free in the complement of a finite subset) cyclic topological actions 
of prime order, the $G$-signature theorem holds most of the key to the
existence. As for smooth finite group actions, a basic question
is what additional restrictions the smooth structures of the $4$-manifolds
may impose on the fixed-point set and the action around it (e.g. 
nonsmoothability of topological actions). 

In this paper, we initiate a study on a class of symmetries of smooth 
$4$-manifolds, which we call {\it symplectic symmetries}. These are smooth
finite group actions which preserve {\it some} symplectic structure on the
$4$-manifold. We recall that a symplectic structure is a closed, 
non-degenerate $2$-form; in particular, the manifold is symplectic.

The study of symplectic structures on smooth $4$-manifolds is one of the
most rapidly growing research area in manifold topology in recent years
(cf. e.g. \cite{Gr, McD, Gom, D2, T3}). One of the fundamental problems
in this study is to understand what restrictions a symplectic structure 
may impose on the underlying smooth structure of the $4$-manifold. Thus
the following seems to be a very natural question: how restrictive the
structure of the fixed-point set and the action around it could be for
symplectic symmetries of a smooth $4$-manifold, and how this may depend
on the underlying smooth structure of the $4$-manifold ? 

Our approach to this question is based on an equivariant version of the
work of C.H. Taubes in \cite{T1, T2} (i.e. the $G$-equivariant 
Seiberg-Witten-Taubes theory). The basic observation is that when the
equivariant versions of Taubes' theorems are applied, the canonical class
of the symplectic $4$-manifold is represented by a set of $2$-dimensional
symplectic subvarieties which is invariant under the group action. In 
principle, the structure of the fixed-point set and the action around it
can be recovered from the induced action on the $2$-dimensional symplectic
subset. This is most effective when the symplectic $4$-manifold is minimal
with $c_1^2=0$; indeed in this case, the set of the $2$-dimensional symplectic
subvarieties representing the canonical class is relatively simple.
While this paper focuses mainly on the structure of the fixed-point set, 
the question as how the symplectic symmetries may depend on the underlying
smooth structure is investigated in \cite{CK}. 

When studying finite group actions on manifolds it is often important and
beneficial to consider an induced action on some algebraic invariants 
associated with the manifold. There is an abundance of such situations, for
example, the notion of a Reidemeister torsion and classification of lens
spaces (cf. \cite{Mil}), the classical Hurwitz theorem about rigidity of 
group actions on surfaces of genus $\geq 2$, rigidity of holomorphic actions
on $K3$ surfaces (cf. \cite{BPV}). The mentioned rigidity of group actions
for the real and complex surfaces asserts that these actions are trivial
if they are homologically trivial, i.e., inducing a trivial action on 
homology. 

In this paper, we shall mainly consider actions which are either 
homologically trivial (over $\Q$ coefficients), or slightly more generally, 
induce a trivial action on the second rational homology. (We remark, 
however, that our method is applicable to much more general situations.)
Our first result is the following rigidity theorem. 

\vspace{3mm}

\noindent{\bf Theorem A:}\hspace{2mm}
{\em Let $M$ be a symplectic $4$-manifold which has trivial canonical class 
{\em(}over $\Q$ coefficients{\em)} and nonzero signature, and obeys 
$b_2^{+}\geq 2$. Then any homologically trivial {\em(}over $\Q$ 
coefficients{\em)}, symplectic action of a finite group $G$ on $M$ must 
be trivial. 
}

\vspace{3mm}

It should be pointed out that if the action is not required to be 
homologically trivial, then there
are nontrivial holomorphic (in particular symplectic) automorphisms of
K\"{a}hlerian $K3$ surfaces (cf. \cite{Nik}). Moreover, the condition
$b_2^{+}\geq 2$ in the theorem is necessary because there exist nontrivial
holomorphic involutions on Enriques surfaces (where $b_2^{+}=1$) which are
homologically trivial (cf. \cite{MN}). Finally, note that the nonvanishing
of signature is also necessary as there is an abundance of homologically
trivial, symplectic finite group actions on the standard symplectic
$4$-torus. 

Note that the $4$-manifold $M$ in Theorem A has a unique Seiberg-Witten basic
class by work of Taubes \cite{T1, T2}. Hence for any other symplectic structure
on $M$ which defines the same orientation, the canonical class is also trivial.
In particular, for the case of $K3$ surfaces (given with the canonical 
orientation), we arrived at the following

\vspace{3mm}

\noindent{\bf Corollary A:} \hspace{2mm}
{\em A homologically trivial symplectic symmetry of a $K3$ surface is
trivial.
}
\vspace{3mm}

\noindent{\bf Remarks:}\hspace{2mm}
There is an open question (due to A. Edmonds, cf. \cite{Kirby}, 
Problem 4.124 (B)) as whether any homologically trivial, smooth actions 
of a finite group on a $K3$ surface must be trivial. The question was 
motivated by the corresponding rigidity for holomorphic actions (cf. 
\cite{BPV}), which was known to be false for topological ones (cf. \cite{E}). 
Corollary A answers the question affirmatively for symplectic symmetries. 
For smooth actions in general, the homological rigidity was only known to 
be true for involutions (cf. \cite{Ma, Rub}). 

\vspace{2mm}

Our main result in this paper is a complete description of the structure
of the fixed-point set for a symplectic cyclic action of prime order on a
minimal symplectic $4$-manifold $M$ with $c_1^2=0$ and $b_2^{+}\geq 2$,
which induces a trivial action on $H^2(M;\Q)$. For simplicity, we shall only 
state the result for the case of pseudofree actions (i.e., actions with only
isolated fixed points). Note that in this case, the induced action on
the tangent space at a fixed point is called the local representation. 
(The general case where the fixed point set may contain $2$-dimensional 
components is addressed in Theorem 3.2, with some refined statements 
contained in Proposition 3.7.)

\vspace{3mm}

\noindent{\bf Theorem B:}\hspace{2mm}{\em
Let $M$ be a minimal symplectic $4$-manifold with $c_1^2=0$ and $b_2^{+}
\geq 2$, which admits a nontrivial, pseudofree action of $G\equiv\Z_p$, 
where $p$ is prime, such that the symplectic structure is preserved under the
action and the induced action on $H^2(M;\Q)$ is trivial. Then the set 
of fixed points of $G$ can be
divided into groups each of which belongs to one of the following five 
possible types. {\em(}We set $\mu_p\equiv\exp(\frac{2\pi i}{p})$.{\em)}
\begin{itemize}
\item [{(1)}] One fixed point with local representation 
$(z_1,z_2)\mapsto (\mu_p^k z_1,\mu_p^{-k}z_2)$ 
for some $k\neq 0\bmod p$, i.e., with local
representation contained in $SL_2(\C)$.
\item [{(2)}] Two fixed points with local representation 
$(z_1,z_2)\mapsto (\mu_p^{2k} z_1,\mu_p^{3k}z_2)$,
$(z_1,z_2)\mapsto (\mu_p^{-k} z_1,\mu_p^{6k}z_2)$ for some $k\neq 0
\bmod p$ respectively. This type of fixed points occurs only when $p>5$.
\item [{(3)}] Three fixed points, one with local representation 
$(z_1,z_2)\mapsto (\mu_p^k z_1,\mu_p^{2k}z_2)$ and the other two with 
local representation $(z_1,z_2)\mapsto (\mu_p^{-k} z_1,\mu_p^{4k}z_2)$ for
some $k\neq 0\bmod p$. This type of fixed points occurs only when $p>3$.
\item [{(4)}] Four fixed points, one with local representation 
$(z_1,z_2)\mapsto (\mu_p^k z_1,\mu_p^k z_2)$ and the other three with 
local representation $(z_1,z_2)\mapsto (\mu_p^{-k} z_1,\mu_p^{3k}z_2)$ for
some $k\neq 0\bmod p$. This type of fixed points occurs only when $p>3$.
\item [{(5)}] Three fixed points, each with local representation 
$(z_1,z_2)\mapsto (\mu_p^k z_1,\mu_p^{k}z_2)$ for some $k\neq 0\bmod p$. 
This type of fixed points occurs only when $p=3$.
\end{itemize}
}

\vspace{2mm}

Combined with the $G$-signature theorem, Theorem B implies the following 
rigidity for the corresponding homologically trivial actions. 

\vspace{3mm}

\noindent{\bf Corollary B:}\hspace{2mm}{\em
Let $M$ be a minimal symplectic $4$-manifold with $c_1^2=0$ and $b_2^{+}
\geq 2$, which admits a homologically trival 
{\em(}over $\Q$ coefficients{\em)}, pseudofree, symplectic $\Z_p$-action
for a prime $p>1$. Then the following conclusions hold.
\begin{itemize}
\item [{(a)}] The action is trivial if $p\neq 1\bmod 4$, $p\neq 1\bmod 6$,
and the signature of $M$ is nonzero. In particular, if the signature of 
$M$ is nonzero, then for infinitely many primes $p$ the manifold
$M$ does not admit any such nontrivial $\Z_p$-actions.
\item [{(b)}] The action is trivial as long as there is a fixed point of 
type {\em(}1{\em)} in Theorem B. 
\end{itemize}
}

\vspace{2mm}

\noindent{\bf Remarks:}\hspace{2mm}
(1) Corollary B shows that symplectic symmetries are far more restrictive 
than topological ones. Indeed, there is the following theorem of Edmonds
(cf. \cite{E}).

\vspace{2mm}

{\it Let $M$ be a closed simply-connected $4$-manifold. Then there is a 
locally linear, pseudofree, homologically trivial action of $\Z_p$ on
$M$ for every prime $p>3$.
}
\vspace{2mm}

On the other hand, one should compare the nonexistence results in 
Corollary B with a theorem of Edmonds (cf. \cite{McC}) which
is of a purely topological nature.

\vspace{2mm}

{\it Let $M$ be a closed simply-connected $4$-manifold with $b_2\geq 3$. 
If a finite group $G$ acts on $M$ locally linearly, pseudofreely, and 
homologically trivially, then $G$ must be a cyclic group.
}

\vspace{2mm}

(2) Corollary B should also be compared with the (stronger) results of Peters 
(cf. \cite{Pe2}) on homological rigidity of holomorphic actions on 
K\"{a}hlerian elliptic surfaces. However, we would like to point out 
that Peters' result relies essentially on the fact that a homologically
trivial holomorphic automorphism on an elliptic surface preserves the 
holomorphic elliptic fibration, which does not have an analog in the 
symplectic category. Moreover, symplectic $4$-manifolds are a much larger 
class of manifolds than K\"{a}hler surfaces (cf. \cite{Gom}) (note that 
this is true even for the case of $c_1^2=0$, cf. \cite{FS1}), and 
symplectic automorphisms are much more flexible than holomorphic 
automorphisms --- the former form an infinite dimensional Lie Group while 
the latter only a finite dimensional one. 

(3) It is essential that the symplectic automorphisms in Theorem A and 
Corollary B are of finite order, i.e. our rigidity for symplectic 
automorphisms is a phenomenon of finite group actions. The necessity of
being finite order can be easily seen from the following example: Let
$\phi:M\rightarrow M$ be any Hamiltonian symplectomorphism. Then the group
$G\equiv \{\phi^n|n\in\Z\}$ acts on $M$ homologically trivially because $\phi$
is homotopic to identity. Of course, $G$ will be of infinite order in general. 

In comparison, the group of holomorphic automorphisms of (compact) 
K\"{a}hler manifolds preserving the K\"{a}hler form is compact 
(cf. \cite{L, Fuj}); in particular, the group is automatically finite 
if the K\"{a}hler manifold admits no nonzero holomorphic 
vectorfields. This is certainly not true in the symplectic category.

(4) Theorem A may be regarded as a special case of Corollary B(b). In
fact, when applying the equivariant versions of Taubes' theorems in 
\cite{T1, T2} in the context of Theorem A, the $2$-dimensional symplectic
subvarieties representing the canonical class must be empty because the
canonical class is trivial in this case. This has the consequence that
the action is pseudofree, and that the fixed-point set consists of only 
type (1) fixed points in Theorem B, which is nonempty because the 
$4$-manifold has nonzero signature.

\vspace{2mm}

The rest of the paper is contained in two sections. The first one, which 
consists of two parts, gives a proof of Theorem A. We present
the proof of Theorem A separately so as to illustrate the 
general philosophy of our paper while keeping the necessary technicalities
at bay. The second section contains discussion of symplectic $\Z_p$-actions
(not necessarily pseudofree) on a minimal symplectic $4$-manifold with
$c_1^2=0$ and $b_2^{+}\geq 2$, and gives a proof for Theorem B and 
Corollary B. It also contains some examples which illustrate our 
considerations.

\vspace{3mm}

\noindent{\bf Acknowledgments:}\hspace{2mm}
We would like to thank Reinhard Schultz for turning our attention to
problems discussed in this paper. We also wish to thank an anonymous referee
whose extensive comments have led to a much improved exposition.

\section{The Proof of Theorem A}

This section is divided into two parts. In part 1 we first briefly 
discuss the $G$-equivariant Seiberg-Witten-Taubes theory, then we use it
to show that the fixed points of the action are isolated, with local
representations all contained in $SL_2(\C)$. In part 2 we combine the
fixed-point data in part 1 with the $G$-signature theorem to show that
the action is trivial.

\vspace{2mm}

\noindent{\bf Part 1:}\hspace{2mm} 
Let $M$ be a smooth, oriented $4$-manifold with an orientation-preserving
action of a finite group $G$. We shall begin by considering the 
$G$-equivariant Seiberg-Witten theory on $M$ (which is equivalent to
the Seiberg-Witten theory on the orbifold $M/G$, cf. \cite{C}). 

More concretely, we fix a Riemannian metric $g$ on $M$ such that $G$
acts via isometries. Suppose there exists a $G$-$Spin^\C$ structure on M,
i.e. a lifting of the pricipal $SO(4)$-bundle of orthonormal frames on
$M$ to a principal $Spin^\C (4)$-bundle as $G$-bundles. Then there are
associated $U(2)$ vector $G$-bundles (of rank $2$) $S_{+}$, $S_{-}$ with
$\det(S_{+})=\det(S_{-})$, and a $G$-equivariant Clifford multiplication
which maps $T^\ast M$ into the skew adjoint endomorphisms of 
$S_{+}\oplus S_{-}$.

With the preceding understood, we consider the $G$-equivariant Seiberg-Witten
equations associated to the $G$-$Spin^\C$ structure 
$$
D_A\psi=0 \mbox{ and } P_{+}F_A=\frac{1}{4}\tau(\psi\otimes\psi^\ast)
+\mu
$$
which are equations for pairs $(A,\psi)$, where $A$ is a $G$-equivariant
$U(1)$-connection on $\det(S_{+})$ and $\psi\in\Gamma(S_{+})$ is a 
$G$-equivariant smooth section of $S_{+}$. As for the notations involved,
here $D_A: \Gamma(S_{+})\rightarrow\Gamma(S_{-})$ is the Dirac operator 
canonically defined from the Levi-Civita connection on $M$ and the
$U(1)$-connection $A$ on $\det(S_{+})$, $P_{+}$ is the orthogonal projection 
onto the subspace of self-dual $G$-equivariant $2$-forms, $\tau$ (which is
canonically defined from the Clifford multiplication) is a map from the space
of endomorphisms of $S_{+}$ into the space of imaginary valued self-dual
$G$-equivariant $2$-forms, and $\mu$ is a fixed, imaginary valued self-dual
$G$-equivariant $2$-form which is added in as a perturbation term. Note that
the $G$-equivariant Seiberg-Witten equations are invariant under the gauge
transformations $(A,\psi)\mapsto (A-2\varphi^{-1}d\varphi,\varphi\psi)$,
where $\varphi$ are circle valued $G$-invariant smooth functions on $M$.

Let $b_G^{2,+}$ be the dimension of the maximal subspace of $H^2(M;\R)$
over which the cup-product is positive and the induced action of $G$ is 
trivial. Then as in the non-equivariant case, the space $\M^G$ of gauge
equivalence classes of the solutions to the $G$-equivariant Seiberg-Witten
equations is compact, and when $b_G^{2,+}\geq 2$, it is an orientable 
smooth manifold (if nonempty) for a generic choice of $(g,\mu)$, whose 
cobordism class is independent of the data $(g,\mu)$. An invariant,
denoted by $SW_G(S_{+})$, can be similarly defined, which, for instance
when $\dim \M^G=0$, is a signed sum of the points in $\M^G$. Moreover,
there is an involution $I$ on the set of $G$-$Spin^\C$ structures which obeys
$\det(I(S_{+}))=-\det(S_{+})$ and $SW_G(I(S_{+}))=\pm SW_G(S_{+})$. 

Given the above terminology and notation, we claim the following

\begin{lemma}
Let $(M,\omega)$ be a symplectic $4$-manifold with trivial canonical class
$c_1(K)$. Suppose a finite group $G$ acts on $M$ via symplectic automorphisms 
such that $b_G^{2,+}\geq 2$. Then the canonical bundle $K$ is isomorphic to 
the trivial bundle $M\times \C$ as $G$-bundles, where $G$ acts on the second 
factor of $M\times \C$ trivially. 
\end{lemma}

\begin{proof}
Fix a $G$-equivariant $\omega$-compatible almost complex structure $J$, 
and let $g=\omega(\cdot,J(\cdot))$ be the associated $G$-equivariant
Riemannian metric. Then there is a canonical $Spin^\C$-structure on $M$
such that the associated $U(2)$ bundles are given by $S_{+}^0=\I\oplus
K^{-1}$ and $S_{-}^0=T^{0,1}M$, where $\I$ is the trivial complex line
bundle. Clearly, this canonical $Spin^\C$-structure is a $G$-$Spin^\C$
structure on $M$, with $\I$ being understood as the $G$-bundle $M\times \C$ 
where $G$ acts on the second factor trivially. 

According to Taubes \cite{T1}, there is a canonical (up
to gauge equivalence) connection $A_0$ on $K^{-1}=\det(S_{+}^0)$, such
that if we set $u_0=(1,0)\in\Gamma(\I\oplus K^{-1})$, then for sufficiently
large $r>0$, $(A_0,\sqrt{r} u_0)$ is the only solution (up to
gauge equivalence) to the Seiberg-Witten equations
$$
D_A\psi=0 \mbox{ and } P_{+}F_A=\frac{1}{4}\tau(\psi\otimes\psi^\ast)
+\mu,
$$
where $\mu=-\frac{i}{4}r\omega+P_{+}F_{A_0}$, and furthermore, 
$(A_0,\sqrt{r} u_0)$ is a non-degenerate solution. In the present situation,
observe that $A_0$ is $G$-equivariant, so is the perturbation 
$\mu=-\frac{i}{4}r\omega+P_{+}F_{A_0}$. Hence $(A_0,\sqrt{r} u_0)$, which
is $G$-equivariant, is also a solution to the $G$-equivariant Seiberg-Witten
equations. We claim that $(A_0,\sqrt{r} u_0)$ is the only solution up to
$G$-invariant gauge equivalence. To see this, suppose $(A,\sqrt{r} u)$
is another solution. Then $(A,\sqrt{r} u)=(A_0-2\varphi^{-1}d\varphi,
\varphi\cdot\sqrt{r} u_0)$ for some circle valued smooth function $\varphi$
on $M$, from which we see that $\varphi$ is $G$-invariant. Hence the claim.
Now note that $(A_0,\sqrt{r} u_0)$ is also non-degenerate as a solution
to the $G$-equivariant Seiberg-Witten equations. This implies that
$SW_G(S_{+}^0)=\pm 1$. 

By the assumption $b_G^{2,+}\geq 2$, we have $SW_G(I(S_{+}^0))=\pm 1$
as well, where $I$ is the involution on the set of $G$-$Spin^\C$
structures, and $I(S_{+}^0)=K\oplus \I$. Thus for any $r>0$, the 
$G$-equivariant Seiberg-Witten equations associated to the $G$-$Spin^\C$
structure $I(S_{+}^0)$, with perturbation 
$\mu=-\frac{i}{4}r\omega+P_{+}F_{A_0}$, has a solution $(A,\psi)$. 
If we write $\psi=\sqrt{r}(\alpha,\beta)\in\Gamma(K\oplus\I)$, then
according to Taubes \cite{T2}, the zero locus $\alpha^{-1}(0)$, if
nonempty, will pointwise converge, as $r\rightarrow +\infty$, to
a set of finitely many $J$-holomorphic curves with multiplicity, which
represents the Poincar\'{e} dual of $c_1(K)$. Since $c_1(K)$ is trivial,
$\alpha^{-1}(0)$ must be empty when $r$ is sufficiently large. The 
section $\alpha$ of $K$ is $G$-equivariant, and is nowhere vanishing for large
$r$, hence it defines isomorphism of $G$-bundles between $K$ and
$M\times \C$, where $G$ acts trivially on the second factor of $M\times\C$.

\end{proof}

Lemma 2.1 provides an important information about the structure of the 
fixed-point set $M^G$ of an action of $G$ on $M$. This information is
summarized in the following

\begin{corollary}
Let $(M,\omega)$ be a symplectic $4$-manifold with trivial canonical class 
$c_1(K)$. Suppose a finite group $G\equiv\Z_p$, $p>1$ prime, acts on $M$ via
symplectic automorphisms such that $b_G^{2,+}\geq 2$. Then the fixed-point
set $M^G\equiv \{m\in M:gm=m, \forall g\in G\}$, if nonempty, consists of
finitely many isolated points, such that with respect to an 
$\omega$-compatible almost complex structure on $M$, the complex
representation of $G\equiv\Z_p$ at each of these fixed points is given
by $(z_1,z_2)\mapsto (\mu_p^k z_1, \mu_p^{-k} z_2)$ for some 
$k\neq 0\bmod p$. {\em(}Here $\mu_p\equiv\exp(\frac{2\pi i}{p})$.{\em)}
\end{corollary}

\begin{proof}
Let $m\in M^G$ be any fixed point. By the equivariant Darboux' theorem,
one can choose an $\omega$-compatible almost complex structure on $M$
which is integrable near $m$, such that there are holomorphic coordinates
$z_1$, $z_2$ near $m$ within which $\omega=\frac{i}{2}(dz_1\wedge d\bar{z}_1
+d z_2\wedge d\bar{z}_2)$ and the action of $G$ is given by  $(z_1,z_2)
\mapsto (\mu_p^{m_1}z_1, \mu_p^{m_2} z_2)$ for some $m_1, m_2$, where 
$m_1\neq 0\bmod p$. Moreover, $dz_1\wedge dz_2$ defines a local section 
of the canonical bundle $K$, with the induced action of $G$ given by
$$
dz_1\wedge dz_2 \mapsto \mu_p^{-(m_1+m_2)} dz_1\wedge dz_2.
$$ 
By Lemma 2.1, $K$ and $M\times \C$ are isomorphic as $G$-bundles, 
where $G$ acts trivially on the second factor of $M\times\C$. 
This implies that $m_1+m_2=0\bmod p$,
from which the corollary follows easily.

\end{proof}

\noindent{\bf Part 2:} For readers' convenience, we first recall a
version of the $G$-signature theorem which will be used in this paper.

Let $M$ be a closed, oriented smooth $4$-manifold, and let $G\equiv\Z_p$
be a cyclic group of prime order $p$ which acts on $M$ effectively via
orientation-preserving diffeomorphisms. Then the fixed point set $M^G$,
if nonempty, will be in general a disjoint union of finitely many isolated 
points and $2$-dimensional orientable submanifolds. Moreover, at each
isolated fixed point, $G$ defines a local complex representation 
$(z_1,z_2)\mapsto (\mu_p^{k}z_1,\mu_p^{kq}z_2)$ for some 
$k, q\neq 0\bmod p$, where $q$ is uniquely determined, and $k$ is 
determined up to a sign. (Here $\mu_p\equiv\exp(\frac{2\pi i}{p})$).

With the preceding understood, we state the $G$-signature theorem 
(cf. \cite{HZ}). 

\begin{theorem}{\em ($G$-signature theorem for prime order cyclic actions)}
$$
|G|\cdot\text{sign}(M/G)=\text{sign}(M)
+\sum_{m}\text{def}_m +\sum_{Y}\text{def}_Y
$$
where $m$ stands for an isolated fixed point, and $Y$ for a 
$2$-dimensional component of $M^G$. The terms, $\text{def}_m$ and 
$\text{def}_Y$ {\em(}which are called  {\em signature defect)}, are given 
by the following formulae:
$$
\text{def}_m=\sum_{k=1}^{p-1}\frac{(1+\mu_p^k)(1+\mu_p^{kq})}
{(1-\mu_p^k)(1-\mu_p^{kq})}
$$
if the local representation at $m$ is given by $(z_1,z_2)\mapsto 
(\mu_p^{k}z_1,\mu_p^{kq}z_2)$, and 
$$
\text{def}_Y=\frac{p^2-1}{3}\cdot (Y\cdot Y)
$$
where $Y\cdot Y$ is the self-intersection of $Y$.
\end{theorem}

Now back to the proof of Theorem A. Without loss of generality we may 
assume that $G\equiv \Z_p$ where $p>1$ is prime. Suppose now that $G$ 
acts trivially on $H_\ast(M;\Q)$. Then $b_G^{2,+}=b_2^{+}\geq 2$ so 
that Corollary 2.2 is true. The Lefschetz fixed point theorem together 
with Corollary 2.2
implies that the fixed-point set $M^G$ is nonmepty and consists of
finitely many isolated points. Indeed, the classical formula relating the
Euler characteristic $\chi(M)$ and the signature $\text{sign}(M)$, i.e.,
$$
2\chi(M)+3\;\text{sign}(M)=c_1^2(K)[M],
$$
gives $\chi(M)=-\frac{3}{2}\;\text{sign}(M)\neq 0$. Moreover, we would like
to point out that the number of fixed points $|M^G|$ equals the Euler
characteristic $\chi(M)$ (cf. \cite{KS}).

On the other hand, by Corollary 2.2 we observe that, in the $G$-signature 
theorem
$$
|G|\cdot\text{sign}(M/G)=\text{sign}(M)+\sum_{m}\text{def}_m
+\sum_{Y}\text{def}_Y,
$$
$\text{def}_m$ is independent of $m$ and is given by
$$
\text{def}_m=\sum_{k=1}^{p-1}\frac{(1+\mu_p^k)(1+\mu_p^{-k})}
{(1-\mu_p^k)(1-\mu_p^{-k})}
$$
and there are no terms $\text{def}_Y$. With $\text{sign}(M/G)=\text{sign}(M)$ 
and $|M^G|=\chi(M)=-\frac{3}{2}\;\text{sign}(M)$, the $G$-signature theorem 
gives rise to
$$
\text{def}_m=\frac{1}{|M^G|}\cdot (p-1)\cdot\text{sign}(M)
=\frac{2}{3}(1-p).
$$
Theorem A follows easily from the following explicit calculation of
$\text{def}_m$, which contradicts the above equation.

\begin{lemma}
$\text{def}_m=\frac{1}{3}(p-1)(p-2)$.
\end{lemma}

\begin{proof} 
It turns out that $\text{def}_m$ can be computed in terms of 
Dedekind sum $s(q,p)$ (cf. \cite{HZ}, page 92), where
$$
s(q,p)=\sum_{k=1}^p((\frac{k}{p}))((\frac{kq}{p}))
$$
with 
$$
((x))=\left\{\begin{array}{ll}
x-[x]-\frac{1}{2} & \mbox{ if } x\in\R\setminus\Z\\
0                 & \mbox{ if } x\in\Z.
\end{array} \right.
$$
(Here $[x]$ stands for the greatest integer less than or equal to $x$.)

In fact, equation (24) in \cite{HZ}, page 180, gives
$$
\text{def}_m=-4p\cdot s(q,p), \mbox{ with } q=-1.
$$
One then computes $6p\cdot s(q,p)$ by (cf. \cite{HZ}, equations (10) and
(9) on page 94)
$$
6p\cdot s(q,p)=(p-1)(2pq-q-\frac{3p}{2})-6f_p(q),
$$
where $f_p(q)=\sum_{k=1}^{p-1} k[\frac{kq}{p}]$. Since $f_p(-1)=
\sum_{k=1}^{p-1} k\cdot (-1)=\frac{1}{2}(1-p)p$, one obtains
$$
\text{def}_m=\frac{1}{3}(p-1)(p-2)
$$
as claimed.
\end{proof}

\begin{remark}

(1) The calculation of $\text{def}_m$ remains valid if $q=-1$ is replaced
by the equivalent (i.e. congruent mod $p$) value of $q=p-1$.

(2) It is quite instructive to compare the proof of homological rigidity
of holomorphic actions on K\"{a}hlerian $K3$ surfaces in \cite{BPV, Pe1} with 
our proof of rigidity of symplectic actions on smooth $4$-manifolds. Both
consists of the same two basic steps. First, one shows that the fixed-point
set of the action consists of isolated fixed points with specific
representations around them. In the holomorphic case this is a rather
obvious observation, whereas in the symplectic case one needs to use 
deep results of Taubes on the Seiberg-Witten theory. Next, the holomorphic
Lefschetz fixed point theorem is used to reach a contradiction for the
$K3$ case. Here the K\"{a}hler condition as well as the holomorphicity of
the action are crucial; in particular the argument uses information about 
the action that is encoded in the quadratic form of the $K3$ surfaces (via
the Hodge decomposition and the holomorphic Lefschetz fixed point theorem).
In our proof the contradiction is reached by computing with the 
$G$-signature theorem. This has the advantage of being applicable to any
smooth (or even locally linear topological) actions on $4$-manifolds.

\end{remark}

\section{Symplectic $\Z_p$-actions on $4$-manifolds with $c_1^2=0$}

In this section we push further the techniques of $G$-equivariant 
Seiberg-Witten-Taubes theory to the case where $c_1^2=0$. The balk of this
section constitutes an analysis of the structure of the fixed-point set
of a symplectic $\Z_p$-action on a minimal symplectic $4$-manifold with
$c_1^2=0$ and $b_2^{+}\geq 2$, which induces a trivial action on the second
rational cohomology. The results are summarized in Theorem 3.2 below (with
some refined statements contained in Proposition 3.7).  
Theorem B is a special case of Theorem 3.2 where the action is pseudofree. 
We set $\mu_p\equiv\exp(\frac{2\pi i}{p})$ throughout. 

Let $(M,\omega)$ be a symplectic $4$-manifold and let $G\equiv\Z_p$, where 
$p>1$ and is prime, act on $M$ via symplectomorphisms. Then the fixed-point
set $M^G$, if nonempty, consists of two types of connected components in 
general: an isolated point or an embedded symplectic surface. Moreover, 
the local representation at each isolated fixed point may be written as 
$(z_1,z_2)\mapsto (\mu_p^{k}z_1,\mu_p^{kq}z_2)$, where the local complex 
coordinates $(z_1,z_2)$ are compatible with the symplectic structure 
$\omega$; in particular, the intergers $k, q$ are uniquely determined in 
their congruence classes (mod $p$). 

Now suppose $b_G^{2,+}\geq 2$ (recall that $b_G^{2,+}$ is the dimension of 
the maximal subspace of $H^2(M;\R)$ over which the cup-product is positive 
and the induced action of $G$ is trivial). Then as we argued in the proof
of Lemma 2.1, given any $G$-equivariant $\omega$-compatible almost complex
structure $J$, there is a solution $(A,\psi)$, with $\psi\in\Gamma(K\oplus
\I)$, to the $G$-equivariant Seiberg-Witten equations with perturbation 
$\mu=-\frac{i}{4}r\omega+P_{+}F_{A_0}$ for any $r>0$. Moreover,
if we write $\psi=\sqrt{r}(\alpha,\beta)$, then as $r\rightarrow +\infty$
the zero set $\alpha^{-1}(0)$ will pointwise converge to a union of finitely
many $J$-holomorphic curves $\cup_i C_i$ such that $c_1(K)$ is Poincar\'{e}
dual to the fundamental class of $\sum_i n_i C_i$ for some $n_i>0$. Since
$\alpha$ is a $G$-equivariant section of $K$, we arrived at the following
observations.
\begin{itemize}
\item The set $\cup_i C_i$ is invariant under the action of $G$.
\item For any fixed point $m\in M^G$, if $m\in M\setminus \cup_i C_i$,
then $\alpha(m)\neq 0$ for sufficiently large $r>0$. Lemma 2.1 and hence
Corollary 2.2 hold true locally near $m$. As a consequence, we see that
$m$ is an isolated fixed point with local representation contained in 
$SL_2(\C)$. 
\end{itemize}

We summarize the preceding discussion in the following

\begin{lemma}
Suppose $b_G^{2,+}\geq 2$. Then given any $G$-equivariant $\omega$-compatible 
almost complex structure $J$, the canonical class $c_1(K)$ 
is represented by the
fundamental class of $\sum_i n_i C_i$ for some $n_i>0$, where $\{C_i\}$
is a finite set of $J$-holomorphic curves which has the following significance.
\begin{itemize}
\item The set $\cup_i C_i$ is invariant under the action of $G$.
\item Let $m$ be a fixed point not contained in $\cup_i C_i$. Then
$m$ must be an isolated fixed point with local representation contained in 
$SL_2(\C)$. In particular, any $2$-dimensional component in $M^G$ is contained
in $\cup_i C_i$.
\end{itemize}
\end{lemma}

Lemma 3.1 allows us to extract information about the fixed-point set $M^G$
and the action around it by analyzing the action of $G$ in a neighborhood
of $\cup_i C_i$. While this requires a priori knowledge about the structure
of $\cup_i C_i$, it can be determined with various additional assumptions
on the canonical class of the manifold. In fact, we will show that when 
$(M,\omega)$ is minimal with $c_1^2=0$, each connected component of 
$\cup_i C_i$ is either a nonsingular elliptic curve, or may be identified 
with a singular fiber of an elliptic fibration. With these understood, 
we now state

\begin{theorem}
Let $(M,\omega)$ be a minimal symplectic $4$-manifold with $c_1^2=0$
and $b_2^{+}\geq 2$, which admits a nontrivial action of $G\equiv\Z_p$, 
where $p$ is prime, such that the symplectic structure is preserved under
the action and the induced action on $H^2(M;\Q)$ is trivial. Then 
for any $G$-equivariant $\omega$-compatible almost complex structure $J$, 
there are $J$-holomorphic curves $\{C_i\}$ and positive integers $\{n_i\}$, 
such that $\cup_i C_i$ is $G$-invariant and the Poincar\'{e} dual of the 
canonical class $c_1(K)$ is represented by the fundamental class of 
$\sum_i n_iC_i$. Furthermore, the following statements describe the 
structure of $\cup_i C_i$ as well as that of the fixed-point set $M^G$.
\begin{itemize}
\item [{(1)}] If a fixed point is not contained in $\cup_i C_i$, 
it must be isolated with local representation contained in $SL_2(\C)$.
\item [{(2)}] The following is a list of all possibilities for a connected 
component of $\cup_i C_i$:
\begin{itemize}
\item [{(I)}] An embedded torus with self-intersection $0$.
\item [{(II)}] A cusp sphere with self-intersection $0$.
\item [{(III)}] A nodal sphere with self-intersection $0$.
\item [{(IV)}] A union of two embedded $(-2)$-spheres intersecting at a 
single point with tangency of order $2$.
\item [{(V)}] A union of embedded $(-2)$-spheres intersecting transversely.
\end{itemize}
\item [{(3)}] Accordingly, the possibilities for the associated fixed-point 
data are listed below, if the connected component of $\cup_i C_i$ contains 
at least one fixed point:
\begin{itemize}
\item [{(i)}] For a type I component, there are three possibilities: (a) 
it is entirely fixed by $G$, (b) it contains four isolated fixed points 
each of which having a local representation contained in $SL_2(\C)$, 
and (c) it contains three isolated fixed points, all having the same local 
representation which is either $(z_1,z_2)\mapsto (\mu_p^k z_1,\mu_p^k z_2)$ 
or $(z_1,z_2)\mapsto (\mu_p^k z_1,\mu_p^{2k} z_2)$ for some $k\neq 0 \bmod p$.
Moreover, case (b) occurs only when $p=2$ and case (c) occurs only when $p=3$.
\item [{(ii)}] For a type II component, there are two isolated fixed points
contained in it. One of them is the cusp-singularity, which has local 
representation $(z_1,z_2)\mapsto (\mu_p^{2k} z_1,\mu_p^{3k} z_2)$ while
the other fixed point has local representation 
$(z_1,z_2)\mapsto (\mu_p^{-k} z_1,\mu_p^{6k} z_2)$, for some $k\neq 0\bmod p$.
This case occurs only when $p\geq 5$.
\item [{(iii)}] A type III component contains only one isolated fixed point, 
which is the nodal point, with local representation contained in $SL_2(\C)$.
\item [{(iv)}] A type IV component contains three isolated fixed points,
one of which is the intersection of the two $(-2)$-spheres. As for local
representations, the intersection point always has 
$(z_1,z_2)\mapsto (\mu_p^k z_1,\mu_p^{2k} z_2)$, and 
for each of the other two fixed points, there are two possibilities: 
$(z_1,z_2)\mapsto (\mu_p^{-k} z_1,\mu_p^{4k} z_2)$ which occurs only when
$p>3$, or $(z_1,z_2)\mapsto (\mu_p^{-k} z_1,\mu_p^{-2k} z_2)$ which occurs 
only when $p=3$, for some $k\neq 0\bmod p$.
\item [{(v)}] For a type V component, there are three possibilities: (a) 
it contains at least one $2$-dimensional component of $M^G$, (b) it contains
only a number of isolated fixed points whose local representations are all 
contained in $SL_2(\C)$, and (c) it contains exactly four isolated fixed 
points, one of which has local representation $(z_1,z_2)\mapsto
(\mu_p^k z_1,\mu_p^k z_2)$ and each of the other three has $(z_1,z_2)
\mapsto (\mu_p^{-k} z_1,\mu_p^{3k} z_2)$ for some $k\neq 0\bmod p$. 
The last possibility occurs only when $p\neq 3$. 
\end{itemize}
\end{itemize}
\end{theorem}

(We would like to point out that the fixed points in Theorem B are slightly
reorganized. For instance, for type (2) fixed points in Theorem B the range
of possible primes $p$ is $p>5$, while in Theorem 3.2 (3) (ii), the range
is $p\geq 5$. This is because when $p=5$, the two fixed points contained in
the type II component all have
local representation in $SL_2(\C)$, therefore are classified in Theorem B
as type (1) fixed points. Similar remarks apply to other cases as well.)

The proof of Theorem 3.2 may be divided into two stages. In the first stage,
we determine the structure of $\cup_i C_i$, while in the second stage, we
analyze the action of $G$ in a neighborhood of $\cup_i C_i$. The following
lemma is the starting point for stage 1. 

\begin{lemma}
Assume $(M,\omega)$ is minimal with $c_1^2=0$. Then $c_1(K)\cdot C_i=0$ 
for all $C_i$, and $C_i^2\geq 0$ unless $C_i$ is an embedded $(-2)$-sphere.  
\end{lemma}

\begin{proof}
Note that for any $C_i$,
$$
c_1(K)\cdot C_i=\sum_j n_j C_j\cdot C_i\geq n_i C_i^2.
$$
Thus if $c_1(K)\cdot C_i<0$, $C_i^2<0$ also, so that
$$
C_1^2+c_1(K)\cdot C_i\leq (-1)+(-1)=-2.
$$
By the adjunction formula (cf. \cite{McD1})
$$
C_i^2 + c_1(K)\cdot C_i+2\geq 2\cdot\text{genus}(C_i),
$$ 
$C_i$ is an embedded sphere with $C_i^2=-1$,
which contradicts the minimality of $(M,\omega)$. Hence $c_1(K)\cdot C_i
\geq 0$ for all $C_i$, which implies $c_1(K)\cdot C_i=0$ for all $C_i$
because $0=c_1^2=\sum_i n_i c_1(K)\cdot C_i$. 

To see the last statement, we note that the adjunction formula with 
$c_1(K)\cdot C_i=0$ gives rise to $C_i^2\geq 2\cdot\text{genus}(C_i)-2$,
which implies that $C_i^2\geq 0$ unless $C_i$ is an embedded $(-2)$-sphere. 

\end{proof}

Let $\{\Lambda_\alpha\}$ be the set of connected components of $\cup_i C_i$,
where we write $\Lambda_\alpha=\cup_{i\in I_\alpha} C_i$ for some index
set $I_\alpha$. Denote by $|I_\alpha|$ the cardinality of $I_\alpha$.

First of all, note that for any $\alpha$ and any $i\in I_\alpha$,
$$
0=c_1(K)\cdot C_i=\sum_{j\in I_\alpha} n_jC_j\cdot C_i=
\sum_{j\neq i} n_jC_j\cdot C_i +n_iC_i^2.
$$
It follows from Lemma 3.3 that if $|I_\alpha|\geq 2$, then $\Lambda_\alpha$ 
is a union of embedded spheres with self-intersection $-2$. To analyze the
case of $|I_\alpha|=1$, we need the following refined version of the 
adjunction formula we used earlier: Let $C$ be a $J$-holomorphic curve, with
$\delta$ double points and a number of branch points indexed by $j$, then
$$
C^2+c_1(K)\cdot C+2=2\cdot\text{genus}(C)+2\delta +\sum_j 2\kappa_j
$$
where $\kappa_j$ denotes the Milnor number of the branch point indexed by $j$
(cf. \cite{MW}, Theorem 7.3). Furthermore, it is well-known that if the 
Milnor number
equals $1$, then the branch point must be the cusp-singularity defined by
the equation $z^2+w^3=0$. Now suppose $|I_\alpha|=1$ and $C$ is the 
$J$-holomorphic curve contained in $\Lambda_\alpha$. Then we have 
$C^2=c_1(K)\cdot C=0$, which implies that $C$ is either an embbeded torus,
or an immersed sphere with one double point (i.e. a nodal sphere),
or a sphere with one cusp-singularity (i.e. a cusp sphere), all with 
self-intersection $0$.

To further analyze the components $\Lambda_\alpha$ with $|I_\alpha|\geq 2$,
we observe the following 

\begin{lemma}
Let $\Lambda_\alpha$ be any connected component such that $|I_\alpha|\geq 2$
and there exist $i,j\in I_\alpha$ with $C_i\cdot C_j\geq 2$. Then 
$|I_\alpha|=2$, and if we denote by $C_1$, $C_2$ the two $J$-holomorphic 
curves contained in $\Lambda_\alpha$, then one of the following is true.
\begin{itemize}
\item [{(1)}] $C_1$, $C_2$ intersect at a single point with tangency of
order $2$.
\item [{(2)}] $C_1$, $C_2$ intersect at two distinct points transversely.
\end{itemize}
\end{lemma}

\begin{proof} 
Note that $0=c_1(K)\cdot C_i=\sum_{k\in I_\alpha} n_k C_k\cdot C_i
\geq n_jC_j\cdot C_i+n_i C_i^2$ implies 
$$
2n_j\leq n_jC_i\cdot C_j\leq -n_i C_i^2=2n_i.
$$
Similarly, one has $2n_i\leq 2n_j$, hence $n_i=n_j$. Moreover, 
one must also have $C_i\cdot C_j=2$ and $|I_\alpha|=2$. The lemma 
follows easily.

\end{proof}

Next we analyze the components $\Lambda_\alpha$ with $|I_\alpha|\geq 2$
and $C_i\cdot C_j=1$ for any distinct $i,j\in I_\alpha$. In this case,
$\Lambda_\alpha$ is a union of embedded $(-2)$-spheres with transverse
intersections. Such a configuration can be conveniently represented by a
graph $\Gamma_\alpha$, where each $C_i\subset \Lambda_\alpha$ corresponds 
to a vertex $v_i\in \Gamma_\alpha$, and each intersection point in 
$C_i\cap C_j$ corresponds to an edge connecting the vertices $v_i$, $v_j$.
Moreover, one can associate a matrix $Q_\alpha=(q_{ij})$, where $q_{ii}=1$
for any $i\in I_\alpha$, $q_{ij}=-\frac{1}{2}$ for any distinct $i,j\in
I_\alpha$ with $C_i\cdot C_j\neq 0$, and $q_{ij}=0$ otherwise. Now observe
that $(\sum_{k\in I_\alpha} n_k C_k)\cdot C_i=c_1(K)\cdot C_i=0$ for all
$i\in I_\alpha$, so that the matrix $Q_\alpha$ satisfies the conditions
(ii), (iii) and (i)' in Lemma 2.10 of \cite{BPV}, which implies that
$Q_\alpha$ is positive semi-definite. By Lemma 2.12 (ii) in \cite{BPV},
the graph $\Gamma_\alpha$ must be one from the list in Figure 1 below
(cf. \cite{BPV}, page 20), with $n\geq 2$ if $\Gamma_\alpha$ is of type
$\tilde{A}_n$.

We shall summarize the analysis on $\cup_i C_i$ by categorizing the 
components $\{\Lambda_\alpha\}$ of $\cup_i C_i$ into the following 
three types.
\begin{itemize}
\item [{(A)}] $|I_\alpha|=1$ and $\Lambda_\alpha$ is either an embbeded 
torus, or a nodal sphere, or a cusp sphere, all with self-intersection $0$.
\item [{(B)}] $|I_\alpha|=2$ and $\Lambda_\alpha$ is a union of two embedded 
$(-2)$-spheres intersecting at a single point with tangency of order $2$.
\item [{(C)}] $|I_\alpha|\geq 2$ and $\Lambda_\alpha$ is a union of 
embedded $(-2)$-spheres intersecting transversely. The corresponding graph
$\Gamma_\alpha$ is one from the list in Figure 1. Note that we allow
$n=1$ in type $\tilde{A}_n$ graphs, which represents case (2) of Lemma 3.4. 
\end{itemize}

\begin{figure}[ht]

\setlength{\unitlength}{1pt}

\begin{picture}(420,420)(-210,-195)

\linethickness{0.5pt}

\put(-170,165){$\tilde{A}_n$}

\put(-120,165){\circle*{3}}

\put(-120,165){\line(1,1){30}}

\put(-90,195){\circle*{3}}

\put(-90,195){\line(1,0){30}}

\put(-60,195){\circle*{3}}

\put(-60,195){\line(1,-1){30}}

\put(-30,165){\circle*{3}}

\put(-120,165){\line(1,-1){30}}

\put(-90,135){\circle*{3}}

\put(-90,135){\line(1,0){30}}

\put(-60,135){\circle*{3}}

\multiput(-60,135)(3,3){10}{\line(1,1){2}}

\put(20,165){($n+1$ vertices, $n\geq 1$)}

\put(-170,90){$\tilde{D}_n$}

\put(-120,105){\circle*{3}}

\put(-120,105){\line(2,-1){30}}

\put(-120,75){\circle*{3}}

\put(-120,75){\line(2,1){30}}

\put(-90,90){\circle*{3}}

\put(-90,90){\line(1,0){30}}

\put(-60,90){\circle*{3}}

\put(-30,90){\circle*{3}}

\multiput(-60,90)(3,0){10}{\line(1,0){2}}

\put(0,90){\circle*{3}}

\put(-30,90){\line(1,0){30}}

\put(0,90){\line(2,1){30}}

\put(0,90){\line(2,-1){30}}

\put(30,105){\circle*{3}}

\put(30,75){\circle*{3}}

\put(60,90){($n+1$ vertices, $n\geq 4$)}

\put(-170,30){$\tilde{E}_6$}

\put(-120,30){\circle*{3}}

\put(-120,30){\line(1,0){30}}

\put(-90,30){\circle*{3}}

\put(-90,30){\line(1,0){30}}

\put(-60,30){\circle*{3}}

\put(-60,30){\line(1,0){30}}

\put(-30,30){\circle*{3}}

\put(-30,30){\line(1,0){30}}

\put(0,30){\circle*{3}}

\put(-60,30){\line(0,-1){60}}

\put(-60,0){\circle*{3}}

\put(-60,-30){\circle*{3}}

\put(-170,-60){$\tilde{E}_7$}

\put(-120,-60){\circle*{3}}

\put(-120,-60){\line(1,0){180}}

\put(-90,-60){\circle*{3}}

\put(-60,-60){\circle*{3}}

\put(-30,-60){\circle*{3}}

\put(0,-60){\circle*{3}}

\put(30,-60){\circle*{3}}

\put(60,-60){\circle*{3}}

\put(-30,-60){\line(0,-1){30}}

\put(-30,-90){\circle*{3}}

\put(-170,-150){$\tilde{E}_8$}

\put(-120,-150){\circle*{3}}

\put(-120,-150){\line(1,0){210}}

\put(-90,-150){\circle*{3}}

\put(-60,-150){\circle*{3}}

\put(-30,-150){\circle*{3}}

\put(0,-150){\circle*{3}}

\put(30,-150){\circle*{3}}

\put(60,-150){\circle*{3}}

\put(90,-150){\circle*{3}}

\put(-60,-150){\line(0,-1){30}}

\put(-60,-180){\circle*{3}}

\end{picture}

\caption{}

\end{figure}

We end the discussion in stage 1 of the proof by identifying the 
$G$-invariant components of $\cup_i C_i$ under the assumption
that $G$ induces a trivial action on $H^2(M;\Q)$.

\begin{lemma}
Assume $G$ induces a trivial action on $H^2(M;\Q)$. Then the following hold
true.
\begin{itemize}
\item [{(1)}] If a component $\Lambda_\alpha$ is not $G$-invariant, then
$\Lambda_\alpha$ must be of type $(A)$.
\item [{(2)}] If a component $\Lambda_\alpha$ is $G$-invariant, then each 
$J$-holomorphic curve $C_i\subset\Lambda_\alpha$ must also be $G$-invariant.
\end{itemize}
\end{lemma}

\begin{proof} 
First, if a component $\Lambda_\alpha$ is not $G$-invariant, then for any
$C_i\subset\Lambda_\alpha$, $C_i$ is disjoint from $g\cdot C_i$ for any
$1\neq g\in G$. Particularly, $C_i^2=(g\cdot C_i)\cdot C_i=0$ because
$G$ induces a trivial action on $H^2(M;\Q)$. Clearly such a component 
$\Lambda_\alpha$ is of type (A).

Secondly, suppose $\Lambda_\alpha$ is $G$-invariant. Then if for some $C_i
\subset\Lambda_\alpha$, $g\cdot C_i\neq C_i$ for some $g\in G$, then 
$|I_\alpha|\geq 2$, and hence $\Lambda_\alpha$ is either of type (B) or
type (C). In any case, $-2=C_i^2=(g\cdot C_i)\cdot C_i\geq 0$, which
is a contradiction.

\end{proof}

Next we enter stage 2 of the proof where we analyze the action of $G$ 
in a neighborhood of $\cup_i C_i$. For this purpose,
we need to introduce the following terminology. Let $C$
be any $G$-invariant $J$-holomorphic curve which is not fixed by $G$.
Suppose $C$ is parametrized by an equivariant $J$-holomorphic map
$f:\Sigma\rightarrow M$, where $\Sigma$ is a Riemann surface with a 
$G\equiv\Z_p$ holomorphic action. Let $z_i\in\Sigma$, $i=1,2,\cdots,k$, 
be the fixed points of $G$, and for each $i$, let $g_i\in G$ be the unique
element whose action near $z_i$ is given by a counterclockwise rotation of
angle $\frac{2\pi}{p}$. Set $p_i=f(z_i)$, which is a fixed point of $G$, 
and let $(m_{i,1},m_{i,2})$, where $0\leq m_{i,1},m_{i,2}<p$, be a pair 
of integers such that the action of $g_i$ on the tangent space at $p_i$
is given by $(w_1,w_2)\mapsto (\mu_p^{m_{i,1}}w_1,\mu_p^{m_{i,2}}w_2)$.
We would like to point out that $(m_{i,1},m_{i,2})$ is uniquely determined 
up to order by the $J$-holomorphic curve $C$, and $f$ is embedded 
near $z_i$ iff one of $m_{i,1}, m_{i,2}$ equals $1$. We shall call
$(m_{i,1},m_{i,2})$ a pair of {\it rotation numbers} at $p_i$ associated 
to $C$. Note that when $C$ has no double points, each fixed point $p_i$ has a
unique pair of rotation numbers associated to $C$.

We shall consider specially the case where $\Sigma$ is a Riemann sphere
and the $J$-holomorphic curve $C$ satisfies $c_1(K)\cdot C=0$. (Note that
every $C_i$ in $\cup_i C_i$ satisfies this condition, cf. Lemma 3.3.)
It is clear that there are exactly two fixed points $z_1, z_2$ of the
$G\equiv \Z_p$ action on $\Sigma$. The following congruence relation
for the rotation numbers will be frequently used later in the proof.

\begin{lemma} Suppose $C$ is a sphere and $c_1(K)\cdot C=0$. Then 
$$
\sum_{i=1}^2(m_{i,1}+m_{i,2})=0\bmod p
$$
where $(m_{i,1},m_{i,2})$, $i=1,2$, are the rotation numbers associated
to $C$.
\end{lemma}

\begin{proof}
To see this, let $f:\s^2\rightarrow M$ be a $G$-equivariant $J$-holomorphic 
parametrization of $C$. Then the virtual dimension of the moduli space of
the corresponding equivariant $J$-holomorphic maps at $f$, which is the 
index of a first order elliptic differential operator of Cauchy-Riemann 
type over the orbifold $\s^2/G$, is of even dimension $2d_f$, where
$$
d_f= -\frac{1}{p}c_1(K)\cdot C +2-\sum_{i=1}^2
\frac{m_{i,1}+m_{i,2}}{p}.
$$
(See the Riemann-Roch theorem for orbit spaces in \cite{AS}, or Lemma 3.2.4
in \cite{CR} for the case of general orbifold Riemann surfaces.) The said
congruence relation on the rotation numbers follows easily from
$d_f\in\Z$ and the assumption $c_1(K)\cdot C=0$. 

\end{proof}

With the preceding understood, the next proposition finishes the proof
of Theorem 3.2, and moreover, it provides a refinement for some of the
statements in Theorem 3.2.

\begin{proposition}
Under the assumptions in Theorem 3.2, the following hold true for any
$G$-equivariant $J$, where $\Lambda_\alpha$ is a connected component of
$\cup_i C_i$ which contains at least one fixed point of $G$.
\begin{itemize}
\item [{(1)}] Let $\Lambda_\alpha=C$ be a $G$-invariant, type $(A)$
component. Then there are three possibilities: 
\begin{itemize}
\item [{(i)}] $C$ is an embedded torus. In this case, $C$ is either 
fixed entirely by $G$, or it contains four isolated fixed points if
$p=2$, or it contains three isolated fixed points if $p=3$. The rotation
numbers at each fixed point are all the same, which is $(1,1)$ if $p=2$,
and either $(1,1)$ or $(1,2)$ if $p=3$.
\item [{(ii)}] $C$ is a cusp sphere. This case occurs only if $p\geq 5$.
There are two isolated fixed points, one is the cusp-singularity and the 
other is a smooth point, with the rotation numbers being $(2,3)$ and 
$(1,p-6)$ respectively if $p>5$, and $(2,3)$ and $(1,4)$ if $p=5$.
\item [{(iii)}] $C$ is a nodal sphere. It contains one isolated fixed 
point, the nodal point, with two pairs of rotation numbers both equaling 
$(1,p-1)$.
\end{itemize}
\item [{(2)}] Let $\Lambda_\alpha$ be a type $(B)$ component. This case 
occurs only when $p\geq 3$. In this case each of the $(-2)$-spheres in 
$\Lambda_\alpha$ contains two isolated fixed points, where one of them 
is the intersection of the two spheres, with the rotation numbers 
associated to either sphere being $(1,2)$, while the other fixed point 
has rotation numbers $(1,p-4)$ if $p>3$, and $(1,2)$ if $p=3$.
\item [{(3)}] Let $\Lambda_\alpha$ be a type $(C)$ component. Then there
are four possibilities:
\begin{itemize}
\item [{(i)}] $\Lambda_\alpha$ contains a $2$-dimensional component of $M^G$, 
with $n=4\bmod p$ if $\Lambda_\alpha$ is represented by a type $\tilde{D}_n$
graph and $n=-1\bmod p$ if $\Lambda_\alpha$ is represented by a type 
$\tilde{A}_n$ graph.
\item [{(ii)}] $\Lambda_\alpha$ is of type $\tilde{A}_n$ and the intersection 
of each pair of spheres is an isolated fixed point, with rotation numbers
$(1,p-1)$ associated to either sphere. 
\item [{(iii)}] $\Lambda_\alpha$ is of type $\tilde{A}_2$ where the three 
spheres intersect at a single point; there are four isolated fixed
points, one occurs at the intersection point and each of the other three
is contained in each one of the three spheres, with the rotation numbers
associated to each sphere being $(1,1)$ at the intersection point and
$(1, |p-3|)$ at each of the other three fixed points. This case occurs 
only if $p\neq 3$.
\item [{(iv)}] $\Lambda_\alpha$ is of type $\tilde{A}_1$ which contains
four isolated fixed points. The rotation numbers 
at each fixed point is $(1,1)$, and this case occurs only if $p=2$.
\end{itemize}
\end{itemize}
\end{proposition}

\begin{proof} 
(1) (i) Suppose $C$ is an embedded torus. Since $C^2=0$, a regular 
neighborhood of $C$ has boundary $T^3$. If $C$ is not fixed by $G$, 
then the induced $G$-action on $T^3$ must be free. Such actions are 
classified (e.g. see Theorem 4.3 and Table 4.4 in \cite{Sco}). In 
particular, since $C$ contains a fixed point by assumption, and $p$ is 
prime, it follows that this happens only if $p=2$ or $p=3$. Moreover, 
the quotient space $T^3/G$ is naturally a Seifert manifold with base 
$S^2$, and with normalized Seifert invariant $(-2, (2,1),(2,1),(2,1),(2,1))$ 
if $p=2$, and $(-1, (3,1),(3,1),(3,1))$ or $(-2, (3,2),(3,2),(3,2))$
if $p=3$. The statement about rotation numbers follows easily from the
description of the normalized Seifert invariant above.

(ii) Suppose $C$ is a cusp sphere. Then $C$ is not fixed by $G$ because
it is not a smooth surface. Moreover, $C$ contains two isolated fixed points,
one of which is the cusp-singularity. Because near the cusp-singularity, $C$
can be parametrized by a $J$-holomorphic map $z\mapsto (z^2,z^3+\cdots)$, it 
is easily seen that the rotation numbers are $(2,3)$. Let $(1,m)$, where
$1<m<p$, be the rotation numbers at the other fixed point, near which $C$
is embedded. Then by Lemma 3.6, the congruence relation $(2+3)+
(1+m)=0\bmod p$ must hold. It follows easily that $p\geq 5$, and moreover,
$m=p-6$ when $p>5$ and $m=4$ when $p=5$. 

(iii) Suppose $C$ is a nodal sphere. 
Again $C$ can not be fixed by $G$ because it 
is immersed, and furthermore, the double point must be a fixed point which 
is easily seen the only fixed point contained in $C$. There are two pairs 
of rotation numbers. A simple inspection shows that they are both $(1,p-1)$. 

(2) Suppose $\Lambda_\alpha$ is a type (B) component. Then by Lemma 3.5,
$\Lambda_\alpha$ is $G$-invariant and so is each of the embedded 
$(-2)$-spheres in $\Lambda_\alpha$. Note that neither of the two 
spheres is fixed by $G$. This is because if one of them is fixed 
by $G$, so is the other as they intersect at a point with tangency of 
order $2$. This, however, is impossible.

It is clear that each of the embedded $(-2)$-spheres contains two isolated 
fixed points where one of the fixed points is the intersection of the two 
spheres. It remains only to check the rotation numbers at each fixed point. 
For the intersection point, if we choose a local complex coordinate system 
$(w_1,w_2)$ such that one of the embedded $(-2)$-spheres is defined by 
$w_2=0$, then the other embedded $(-2)$-sphere is locally parametrized 
by $z\mapsto (z,z^2+\cdots)$. It follows easily that the rotation numbers 
associated to either embedded $(-2)$-sphere is $(1,2)$ at the intersection 
point. On the other hand, if we let $(1,m)$ be the rotation numbers at the 
other fixed point on the sphere, then the congruence relation in Lemma 3.6
gives $(1+2)+(1+m)=0\bmod p$, which implies that $p\geq 3$, and moreover, 
$m=p-4$ if $p>3$ and $m=2$ if $p=3$.

(3) Suppose $\Lambda_\alpha$ is a type (C) component. Then by Lemma 3.5,
each embedded $(-2)$-spheres in $\Lambda_\alpha$ is $G$-invariant. There
are two possible scenarios: (a) a regular neighborhood of $\Lambda_\alpha$
in $M$ is a plumbing of embedded $(-2)$-spheres, (b) $\Lambda_\alpha$
consists of three embedded $(-2)$-spheres which intersect at a single point.

Let's first consider scenario (a). It is clear that if there are two distinct
spheres in $\Lambda_\alpha$ with one intersection point not fixed by $G$,
then we must have $p=2$ and $\Lambda_\alpha$ is of type $\tilde{A}_1$, with
each sphere containing two isolated fixed points of rotation numbers $(1,1)$.
(This case is listed as (3) (iv) in the proposition.)
Now suppose the intersection of any two distinct spheres in $\Lambda_\alpha$
is fixed by $G$. Then a sphere in $\Lambda_\alpha$ which intersects with more 
than two other spheres must be entirely fixed by $G$, because a $\Z_p$-action
on $\s^2$ can not have more than two fixed points unless it is trivial.
In particular, $\Lambda_\alpha$ contains a $2$-dimensional component of
$M^G$ if the graph representing $\Lambda_\alpha$ is of type 
$\tilde{D}_n$, $\tilde{E}_6$, $\tilde{E}_7$ or $\tilde{E}_8$.

We need to further discuss the cases where $\Lambda_\alpha$ is of type
$\tilde{A}_n$ or $\tilde{D}_n$. 

First we assume $\Lambda_\alpha$ is of type $\tilde{D}_n$, and show that
$n=4\bmod p$. The case where $n=4$ is trivial, so we assume $n>4$. Then 
there are two vertices in the representing graph $\Gamma_\alpha$, each of
which connects to three other vertices, and moreover, there is a linear
sub-graph consisting of $(n-3)$ vertices, with these two vertices at each
end. We denote by $v_0,v_1,\cdots,v_k$ the vertices along the linear 
sub-graph, where $k=n-4$.  

According to \cite{Or}, a regular neighborhood of the configuration of
embedded $(-2)$-spheres represented by the linear sub-graph may be obtained 
by an $\s^1$-equivariant plumbing where the $\s^1$-action is linear on each
sphere. On the other hand, the action of $G$ on each embedded $(-2)$-sphere 
is conjugate to a linear action. It follows that there is an 
$\s^1$-equivariant plumbing such that 
the $G\equiv\Z_p$ action on the regular neighborhood is induced by the 
inclusion $\Z_p\subset \s^1$. The congruence relation $n=4\bmod p$ follows 
from this consideration. 

More concretely, the $\s^1$-equivariant plumbing is done as follows. 
For each $i$, $0\leq i\leq k$, regard the normal bundle of the embedded 
$(-2)$-sphere represented by the vertex $v_i$ as the result of 
$\s^1$-equivariantly sewing $D^2\times D^2$ to $D^2\times D^2$ (here
the second $D^2$ in each $D^2\times D^2$ represents the fiber of the
normal bundle) by the matrix
$$
\left (\begin{array}{ll}
-1 & 0\\
2 & 1\\
\end{array} \right ),
$$
where in polar coordinates on each factor $D^2$, the $\s^1$-action
on the first $D^2\times D^2$ is given by $(r,\gamma,s,\delta)\mapsto 
(r,\gamma+u_{i,1}\theta, s,\delta+v_{i,1}\theta)$, $0\leq \theta\leq 2\pi$,
and on the second $D^2\times D^2$ it is given by $(r,\gamma,s,\delta)\mapsto 
(r,\gamma+u_{i,2}\theta, s,\delta+v_{i,2}\theta)$, $0\leq \theta\leq 2\pi$,
for some $u_{i,1}$, $v_{i,1}$, $u_{i,2}$, $v_{i,2}\in\Z$ with
$u_{0,1}=0$, $v_{0,1}=1$. Then the plumbing identifies the second
$D^2\times D^2$ associated to the vertex $v_i$ to the first $D^2\times D^2$
associated to the next vertex $v_{i+1}$ with the two factors of $D^2$ 
switched. Moreover, the plumbing is equivariant so that $u_{i,2}=v_{i+1,1}$,
$v_{i,2}=u_{i+1,1}$ must hold.

With the preceding understood, the weights of the $\s^1$-action on each
$D^2\times D^2$ can be determined from the following equation where $i\geq 0$
(cf. \cite{Or}, page 27)
$$
\left (\begin{array}{l}
u_{i,2}\\
v_{i,2}
\end{array}\right )
=
\left (\begin{array}{ll}
-(i+1) & 1-(i+1)\\
i+2 & i+1\\
\end{array}\right )
\left (\begin{array}{l}
u_{0,1}\\
v_{0,1}
\end{array}\right ).
$$
On the other hand, the $G\equiv\Z_p$ action on $\Lambda_\alpha$ is induced
by $\Z_p\subset\s^1$, and the embedded $(-2)$-sphere represented by $v_k$
is fixed by $G\equiv \Z_p$, so that $u_{k,2}=0\bmod p$ must hold. With
$u_{0,1}=0$ and $v_{0,1}=1$, one has $k=-u_{k,2}=0\bmod p$, or
equivalently, $n=4\bmod p$. This finishes the proof for the case when
$\Lambda_\alpha$ is of type $\tilde{D}_n$.

Next we consider the case where $\Lambda_\alpha$ is of type $\tilde{A}_n$.
We will first show that $n$ satisfies the congruence relation $n=-1\bmod p$
(i.e. the number of vertices in $\Gamma_\alpha$ is divisible by $p$) if one
of the following conditions are satisfied: (1) one of the embedded 
$(-2)$-spheres in $\Lambda_\alpha$ is fixed by $G$, or (2) one of the
isolated fixed points in $\Lambda_\alpha$ has a pair of rotation numbers
$(m_1,m_2)$ such that $m_1+m_2\neq 0\bmod p$. 

Consider the linear graph $\Gamma_\alpha^\prime$
obtained from $\Gamma_\alpha$ by removing any one of its edges. Denote by
$v_1,v_2,\cdots,v_k$ the vertices along $\Gamma_\alpha^\prime$, where
$k=n+1$. Then $n=-1\bmod p$ is equivalent to $k=0\bmod p$, which is what 
we will show next. To see this, note that as we argued earlier, there exists 
an $\s^1$-equivariant plumbing associated to $\Gamma_\alpha^\prime$, such 
that the $G\equiv\Z_p$ action is induced by the 
inclusion $\Z_p\subset\s^1$ from the $\s^1$-action associated to the 
equivariant plumbing. Now suppose at the vertex $v_i$ the $\s^1$-actions 
on the two copies of $D^2\times D^2$ are given (in polar coordinates on each
$D^2$ factor) by $(r,\gamma,s,\delta)\mapsto 
(r,\gamma+u_{i,1}\theta, s,\delta+v_{i,1}\theta)$, $0\leq \theta\leq 2\pi$,
and  $(r,\gamma,s,\delta)\mapsto 
(r,\gamma+u_{i,2}\theta, s,\delta+v_{i,2}\theta)$, $0\leq \theta\leq 2\pi$,
for some integers $u_{i,1}$, $v_{i,1}$, $u_{i,2}$, $v_{i,2}$. Then
$$
\left (\begin{array}{l}
u_{k,2}\\
v_{k,2}
\end{array}\right )
=
\left (\begin{array}{ll}
-k & 1-k\\
k+1 & k\\
\end{array}\right )
\left (\begin{array}{l}
u_{1,1}\\
v_{1,1}
\end{array}\right )
$$
must be satisfied (cf. \cite{Or}, page 27). On the other hand, because the 
$G\equiv\Z_p$ action on $\Lambda_\alpha$ is induced by the inclusion
$\Z_p\subset\s^1$ from the $\s^1$-action associated to the equivariant 
plumbing, it follows that the following equation must hold in 
congruence mod $p$:
$$
\left (\begin{array}{l}
u_{1,1}\\
v_{1,1}\\
\end{array}\right )
=
\left (\begin{array}{ll}
0 & 1\\
1 & 0\\
\end{array}\right )
\left (\begin{array}{l}
u_{k,2}\\
v_{k,2}\\
\end{array}\right )
=
\left (\begin{array}{ll}
k+1 & k\\
-k & 1-k\\
\end{array}\right )
\left (\begin{array}{l}
u_{1,1}\\
v_{1,1}
\end{array}\right )
$$
which gives rise to the congruence relation $k(u_{1,1}+v_{1,1})=0 \bmod p$. 
Hence $k=0\bmod p$ as long as $u_{1,1}+v_{1,1}\neq 0\bmod p$.

To see that $u_{1,1}+v_{1,1}\neq 0\bmod p$, note that if one of the embedded
$(-2)$-spheres in $\Lambda_\alpha$ is fixed by $G$, which without loss of
generality may be assumed to be the one represented by the vertex $v_1$,
then $u_{1,1}=0$ and $v_{1,1}=1$, so that $u_{1,1}+v_{1,1}\neq 0\bmod p$.
If one of the isolated fixed points in $\Lambda_\alpha$ has a pair of
rotation numbers $(m_1,m_2)$ with $m_1+m_2\neq 0\bmod p$, then one may 
similarly assume that $u_{1,1}=m_1$, $v_{1,1}=m_2$, which also implies
$u_{1,1}+v_{1,1}\neq 0\bmod p$.

Therefore for scenario (a) it remains to show that if none of the embedded 
$(-2)$-spheres in $\Lambda_\alpha$ is fixed by $G$, then the rotation numbers
(associated to either sphere) are $(1,p-1)$ at each fixed point. In
particular, we will rule out the second possibility in the preceding
discussion, i.e., one of the isolated fixed points in $\Lambda_\alpha$ has 
a pair of rotation numbers $(m_1,m_2)$ such that $m_1+m_2\neq 0\bmod p$. 

We begin by introducing the following notations. Let $p_i$, $1\leq i\leq k-1$, 
be the intersection of the $i$-th sphere (i.e. the one represented by the 
vertex $v_i$) with the $(i+1)$-th sphere, and $p_k$ be the intersection of the
$k$-th sphere with the first sphere. Let $(1,m_i)$, $(1,m_i^\prime)$
be the rotation numbers associated to the $i$-th sphere at $p_i$, $p_{i-1}$
respectively (here $p_0=p_k$). 

First of all, some basic properties of the integers $m_i$, $m_i^\prime$,
$1\leq i\leq k$. Note first that $m_i^\prime$ and $m_{i-1}$ are mutually
determined by each other in the congruence equation $m_i^\prime m_{i-1}
=1\bmod p$. Secondly, for each $i$, $m_i^\prime$ and $m_i$ satisfy the
congruence relation $(1+m_i^\prime)+(1+m_i)=0\bmod p$ (cf. Lemma 3.6), 
which implies that either $m_i=m_i^\prime=p-1$ or $2+m_i^\prime+m_i=p$. 
Finally, with $m_i^\prime m_{i-1}=1\bmod p$, we see easily that exactly 
one of the following is true:
\begin{itemize}
\item $m_i=m_i^\prime=p-1$ for all $1\leq i\leq k$, or
\item $2+m_i^\prime+m_i=p$ for all $1\leq i\leq k$.
\end{itemize}

Next, we show that the second case, i.e., $2+m_i^\prime+m_i=p$ for all 
$1\leq i\leq k$, can not occur. To see this, we first observe that the
sequence $m_1,m_2, \cdots, m_k$ is periodic with a period $l\leq p-3$.
This is because (1) $1+m_i\neq 0\bmod p$ implies that $k=0\bmod p$, so
that $p\leq k$, (2) each $m_i$ satisfies $1\leq m_i\leq p-3$, hence there
exist $i$, $l$ with $1\leq i\leq k$, $1\leq l\leq p-3$ such that 
$m_i=m_{i+l}$, (3) if $m_i=m_{i+l}$ for some $i$ and $l$, then it holds
for all $i$ with that same $l$. Secondly, a contradiction is reached by
showing that the period $l=1$. Indeed, $l=1$ means that 
$m_1=\cdots =m_k=m$ for some $m$ with $1\leq m\leq p-3$. If we let
$m^\prime$ be the unique integer such that $m^\prime m=1\bmod p$ and
$1<m^\prime <p$. Then $m^\prime+m+2=p$ holds, which implies $m=p-1$,
a contradiction. 

It thus remains to show that $l=1$. To see this, recall that if we do
$\s^1$-equivariant plumbing on the linear graph $\Gamma_\alpha^\prime$
of vertices $v_1,v_2,\cdots,v_k$ with $u_{1,1}=1$ and $v_{1,1}=m_1^\prime$,
then for each $i\leq k$, $u_{i,2}$, $v_{i,2}$ are related to $u_{1,1}$,
$v_{1,1}$ by (cf. \cite{Or}, page 27)
$$
\left (\begin{array}{l}
u_{i,2}\\
v_{i,2}
\end{array}\right )
=
\left (\begin{array}{ll}
-i & 1-i\\
i+1 & i\\
\end{array}\right )
\left (\begin{array}{l}
u_{1,1}\\
v_{1,1}
\end{array}\right ).
$$
On the other hand, we have $v_{i,2}=u_{i,2} m_i\bmod p$, because the
$G\equiv\Z_p$ action on $\Lambda_\alpha$ is the restriction of the
$\s^1$-action associated to the equivariant plumbing to the subgroup
$\Z_p\subset\s^1$. By taking $i=p$ (recall that $p\leq k$), we obtain 
$u_{p,2}=v_{1,1}\bmod p$ and $u_{p,2}m_p=v_{p,2}=u_{1,1}\bmod p$. With
$u_{1,1}=1$, $v_{1,1}=m_1^\prime$, we see that $m_1^\prime m_p=1\bmod p$,
which implies that $m_{p+1}^\prime=m_1^\prime$, and hence $m_{p+1}=m_1$.
This last equality shows that the period $l$ is a divisor of $p$. Hence
$l=1$ because $l\leq p-3$ and $p$ is prime.

This shows that either one of the embedded $(-2)$-spheres in $\Lambda_\alpha$
is fixed by $G$, or the rotation numbers are $(1,p-1)$ at each of the 
isolated fixed points in $\Lambda_\alpha$. 

To complete the proof of Proposition 3.7, it remains to consider scenario (b)
where $\Lambda_\alpha$ consists of three spheres intersecting at a single
point. In this case, it is clear that the intersection point must be fixed 
by $G$. Moreover, since the induced action of $G$ on the tangent space at
this point has three distinct eigenspaces, it must be an isolated
fixed point and has rotation numbers $(1,1)$ associated to each sphere. 
There are three other isolated fixed points in $\Lambda_\alpha$, with each
sphere containing one of them. The rotation numbers are $(1,1)$ if
$p=2$,  and $(1,p-3)$ if $p\neq 2$ by the congruence relation in Lemma 3.6.
Note that scenario (b) does not occur if $p=3$, because otherwise the
rotation numbers $(1,p-3)$ would become $(1,0)$, implying that the normal
direction is fixed by $G$ at the isolated fixed point. 

\end{proof}

We have thus proved Theorem 3.2, and particularly Theorem B in the
introduction, which is a special case where the action is pseudofree.

\vspace{3mm}

With the fixed-point set data in hand, we next give a proof of Corollary B
by appealing to the $G$-signature theorem (cf. Theorem 2.3).

\vspace{3mm}

\noindent{\bf Proof of Corollary B}

\vspace{3mm}

First of all, because $c_1^2=0$ and $G$ acts trivially on $H_\ast (M;\Q)$,
one has
$$
\text{sign}(M)=-\frac{2}{3}\chi(M)=-\frac{2}{3}|M^G|,
$$
where $|M^G|$ is the number of fixed points. Consequently, the $G$-signature
theorem may be rewritten as
$$
-\frac{2}{3}(p-1)\cdot |M^G|=-\frac{2}{3}(p-1)\cdot \chi(M)=
\sum_{m\in M^G} \text{def}_m.
$$

We shall discuss separately two cases: (1) $p=2$ or $p=3$, and (2) $p>3$.

Consider first the case $p=2$ or $p=3$. We remark that in this case we do
not need to use Theorem B, and the corresponding rigidity for the 
$\Z_p$-action is even true for locally linear topological actions,
cf. e.g. \cite{E}. For the sake of completeness, we give a proof here.

If $p=2$, the local representation at each $m\in M^G$ is of the same type,
which is $(z_1,z_2)\mapsto (\mu_p^k z_1,\mu_p^{-k} z_2)$ with $k\neq 0\bmod p$.
Hence $\text{def}_m=\frac{1}{3}(p-1)(p-2)=0$ by Lemma 2.4. This contradicts
the $G$-signature theorem when $\text{sign}(M)\neq 0$. The case where
$p=3$ is similar. There is another type of local representation
$(z_1,z_2)\mapsto (\mu_p^k z_1,\mu_p^k z_2)$ with $k\neq 0\bmod p$, for
which a similar calculation as in Lemma 2.4 shows that 
$\text{def}_m =-\frac{1}{3}(p-1)(p-2)$. For both types of local
representations, $\text{def}_m>-\frac{2}{3}(p-1)$, which contradicts the
$G$-signature theorem when $|M^G|=-\frac{2}{3}\;\text{sign}(M)\neq 0$.

For the rest of the proof, we assume $p>3$. Then the set of fixed points
is divided into groups of the first four types in Theorem B. We introduce
the following notation. For $1\leq k\leq 4$, let $\delta_k$ be the number
of groups of type (k) in Theorem B, and let $\text{def}_{(k)}$ be the
total signature defect contributed by one group of type (k) (i.e., the
sum of $\text{def}_m$ with $m$ running over one group of type (k) fixed 
points). With this notation, the $G$-signature theorem maybe be written as
$$
-\frac{2}{3}(p-1)(\delta_1+2\delta_2+3\delta_3+4\delta_4)=
\delta_1\cdot \text{def}_{(1)}+\delta_2\cdot \text{def}_{(2)}+
\delta_3\cdot \text{def}_{(3)}+\delta_4\cdot \text{def}_{(4)}.
$$
A contradiction will be reached if for all $k=1,2,3$ and $4$,
$$
\text{def}_{(k)}\geq -\frac{2k}{3}(p-1),
$$
with the strict inequality holding for some $k$ with $\delta_k>0$.

The following lemma gives an explicit formula for $\text{def}_{(k)}$,
$k\leq 4$.

\begin{lemma}
\begin{itemize}
\item [{(1)}] $\text{def}_{(1)}=\frac{1}{3}(p-1)(p-2)$ for all $p>1$.
\item [{(2)}] $\text{def}_{(2)}=-8r$ if $p=6r+1$, $\text{def}_{(2)}=8r+8$
if $p=6r+5$.
\item [{(3)}] $\text{def}_{(3)}=-8r$ if $p=4r+1$, $\text{def}_{(3)}=2$
if $p=4r+3$.
\item [{(4)}] $\text{def}_{(4)}=-8r$ if $p=3r+1$, $\text{def}_{(4)}=-4r$
if $p=3r+2$.
\end{itemize}
\end{lemma}

The proof of Lemma 3.8 is given in the appendix, which is a result of
direct calculation. Accepting Lemma 3.8
momentarily, we shall complete the proof of Corollary B for the case
$p>3$. Observe that $\text{def}_{(1)}>-\frac{2}{3}(p-1)$ for all
$p>1$, and $\text{def}_{(2)}\geq -\frac{4}{3}(p-1)$ with equality
only if $p=1\bmod 6$, $\text{def}_{(3)}\geq -\frac{6}{3}(p-1)$ with
equality only if $p=1\bmod 4$, and $\text{def}_{(4)}\geq -\frac{8}{3}(p-1)$ 
with equality only if $p=1\bmod 6$ (note that if $p=3r+1$ and $p$ is
prime, then $r$ must be even and $p=1\bmod 6$). Part (b) of Corollary B
follows from this immediately. As for part (a), observe that
$\text{def}_{(k)}>-\frac{2k}{3}(p-1)$ for all $k$ if $p\neq 1\bmod 4$,
$p\neq 1\bmod 6$. Hence when $\text{sign }(M)\neq 0$, one of $\delta_k$
is nonzero, from which part (a) follows.

\hfill $\Box$

We close this section with two examples. The first one shows that the 
``more exotic'' types of local representations in Theorem B indeed can
occur, at least topologically.

\begin{example}
Let $M$ be a homotopy $K3$ surface (i.e., a manifold homotopy 
equivalent to a $K3$ surface), which is given with the canonical orientation 
such that $\text{sign}(M)=-16$. In this example, we will show that there
are locally linear, homologically trivial topological actions of $\Z_5$
and $\Z_7$ on $M$ with the fixed-point set consisting entirely of type (3) and
type (2) or (4) fixed points in Theorem B respectively. (Compare Corollary 
B(a).) Our construction is based on the work of Edmonds and Ewing \cite{EE}
concerning realization of certain fixed-point data by a locally linear,
topological $\Z_p$-action of prime order on a simply-connected $4$-manifold.

Note that a homologically trivial pseudofree action of $\Z_p$
on $M$ must have $24$ fixed points. Consider first the following fixed-point 
data where $p=5$: Pick $24$ points of $M$, divide them evenly into two 
groups, and assign the points in each group with local representations
$$
(z_1,z_2)\mapsto (\mu_p^k z_1,\mu_p^{2k}z_2),\hspace{2mm}
(z_1,z_2)\mapsto (\mu_p^{-k} z_1,\mu_p^{4k}z_2), \mbox{ and }
(z_1,z_2)\mapsto (\mu_p^{-k} z_1,\mu_p^{4k}z_2)
$$
evaluated at $k=1,2,3,4$.

In order to realize the above fixed-point data, we recall the GSF condition 
from \cite{EE}, which in the present case becomes 
$$
2\;\text{def}_{(3)}=\text{sign}(g,M), \hspace{2mm} \forall g\in\Z_5.
$$
To verify the GSF condition, we note that the action is assumed to be 
homologically trivial, so that $\text{sign}(g,M)=\text{sign}(M)=-16$
for any $g\in\Z_5$. On the other hand, we have $\text{def}_{(3)}=-8$
by Lemma 3.8. Hence the GSF condition is satisfied. Now for $p=5$,
GSF is the only condition needed for the realization of the fixed-point
data by a homologically trivial action, cf. \cite{EE}, Corollary 3.2. 
Consequently, there is a locally linear topological action of $\Z_5$
on $M$ with the fixed-point set consisting entirely of type (3) fixed 
points in Theorem B.

Similar arguments lead to a locally linear, homologically trivial topological 
action of $\Z_7$ on $M$ with the fixed-point set consisting entirely of 
type (2) or type (4) fixed points.

\end{example}

The purpose of the second example is to illustrate that when a certain
additional information about the canonical class and the symplectic
structure is available, Theorem 3.2 (together with Proposition 3.7) may 
yield rigidity for actions which are not necessarily pseudofree.

\begin{example}
Consider a symplectic $4$-manifold $(M,\omega)$ with the following
properties: $M$ is a homotopy $K3$ surface, $\omega$ defines an
integral class $[\omega]\in H^2(M;\R)$ (or more generally, any 
sufficiently small perturbation of an integral class), $c_1^2=0$, and 
$c_1(K)\cdot [\omega]<7$. Note that $(M,\omega)$ must be minimal 
because $M$ has an even intersection form. Moreover, $c_1^2=0$ implies
$\text{sign}(M)=-16$, so that $b_2^{+}=3$.  

\vspace{1.5mm}

We shall next prove that:

\vspace{1.5mm}

{\em There are no nontrivial homologically trivial actions {\em(}not
necessarily pseudofree{\em)} of a finite group on $M$ which preserve the
symplectic structure $\omega$.
}

\vspace{1.5mm}

To see this, suppose there is such an action. Without loss of generality,
we may assume that the action is cyclic
of prime order $p$. By Theorem 3.2 and Proposition 3.7, there is a set of
$J$-holomorphic curves $\{C_i\}$ as described therein, such that $c_1(K)$
is represented by the fundamental class of $\sum_i n_i C_i$ for some
integers $n_i\geq 1$. Since $[\omega]$ is integral, $[\omega]\cdot C_i\geq 1$
for all $C_i$. Now the property $c_1(K)\cdot [\omega]<7$ implies that
in $\cup_i C_i$ there is no component $\Lambda_\alpha$ which is a union
of embedded $(-2)$-spheres represented by a graph of type other than 
$\tilde{A}_n$ or $\tilde{D}_4$. With this understood, it follows easily that
there are only two types of $2$-dimensional components in the fixed-point
set, (1) an embedded torus, (2) an embedded $(-2)$-sphere contained in
a component $\Lambda_\alpha$ which is represented either by a type
$\tilde{A}_n$ graph or by a type $\tilde{D}_4$ graph. Notice further
that a toroidal fixed component does not make any contribution to the
signature defect because it has self-intersection $0$ (cf. Theorem 2.3).

With the preceding understood, we discuss separately the cases (i) the action
is pseudofree, (ii) the action has a $(-2)$-sphere fixed component. 

For case (i), first note that we may assume that there are no fixed points 
whose local representation is contained in $SL_2(\C)$, because otherwise the 
action is trivial by Corollary B(b). Then we see from Theorem 3.2 and 
Proposition 3.7 that the number of fixed points, which equals $\chi(M)$ since 
the action is homologically trivial, equals $2\delta_2+3\delta_3+4\delta_4$, 
where $\delta_2$ is the number of cusp-sphere components $\Lambda_\alpha$, 
$\delta_3$ is the number of $\Lambda_\alpha$ which is a union of two 
$(-2)$-spheres, and $\delta_4$ is the number of $\Lambda_\alpha$ which 
is a union of three spheres. A contradiction is reached easily by observing 
that $\chi(M)=24$, and that $\delta_2+2\delta_3+3\delta_4\leq c_1(K)\cdot 
[\omega]<7$.

For case (ii), we first assume that there is a component $\Lambda_\alpha$
which contains a fixed $(-2)$-sphere and is represented by a type 
$\tilde{A}_n$ graph. Then there will be no component represented by
a type $\tilde{D}_4$ graph because of the constraint 
$c_1(K)\cdot [\omega]<7$, and moreover, by Proposition 3.7 (3), the 
number of $(-2)$-spheres in such a component $\Lambda_\alpha$ is divisible
by the order $p$ of the action. It follows from 
$c_1(K)\cdot [\omega]<7$ that $p=2,3$ or $5$.

Note that in the context of this example, the $G$-signature theorem is 
equivalent to the following equation
$$
-16 (p-1)=\sum_{m}\text{def}_m +\sum_{Y}\text{def}_Y,
$$
where $m, Y$ represent an isolated fixed point and a fixed $(-2)$-sphere
respectively. 

If $p=2$, then $\text{def}_Y=\frac{p^2-1}{3}\cdot (Y\cdot Y)=-2$ for any $Y$
and the definition of $\text{def}_m$ implies that $\text{def}_m=0$ for all 
$m$. We reach a contradiction to the $G$-signature theorem because there 
are at most $3$ fixed $(-2)$-spheres.

Suppose $p=3$ next. Let $\delta_m, \delta_Y$ be the number of isolated fixed
points $m$ and fixed $(-2)$-spheres $Y$ respectively. Then $\delta_m, 
\delta_Y$ obey $\delta_m+2\delta_Y=\chi(M)=24$ (cf. \cite{KS}), and 
moreover, it is clear that $\delta_Y\leq 2$. As for the signature defects
$\text{def}_m$ and $\text{def}_Y$, we note that for $p=3$ there are two
types of local representations: $(z_1,z_2)\mapsto (\mu_p^k z_1,
\mu_p^{-k} z_2)$ or $(z_1,z_2)\mapsto (\mu_p^k z_1,\mu_p^{k} z_2)$. 
A similar calculation as in Lemma 2.4 shows that in both cases, one has
$
\text{def}_m\geq -\frac{1}{3}(p-1)(p-2)=-\frac{2}{3}.
$
On the other hand, $\text{def}_Y=\frac{p^2-1}{3}\cdot (Y\cdot Y)
=-\frac{16}{3}$, which gives a contradiction to the $G$-signature theorem
$$
-32=-16 (p-1)\geq \delta_m\cdot (-\frac{2}{3})+\delta_Y\cdot (-\frac{16}{3})
\geq 24\cdot (-\frac{2}{3})+2\cdot (-\frac{16}{3})=-\frac{80}{3}.
$$

Finally, let $p=5$. First, we observe that there is exactly one component
$\Lambda_\alpha$ of type $\tilde{A}_4$ which contains a fixed $(-2)$-sphere.
Next we recall that if a $(-2)$-sphere is not fixed, then it must contain
exactly two fixed points, and if we let $(1,m_1)$, $(1,m_2)$ be the
corresponding rotation numbers, where $0\leq m_1, m_2<p=5$, the congruence
relation $(1+m_1)+(1+m_2)=0\bmod 5$ must be satisfied (cf. Lemma 3.6). 
A simple inspection shows that $\Lambda_\alpha$ contains, besides the 
fixed $(-2)$-sphere, three isolated fixed points of local representation 
$(z_1,z_2)\mapsto (\mu_p^k z_1,\mu_p^{kq} z_2)$ for some $k\neq 0\bmod p$, 
with $q=1,2,3$ respectively. On the other hand, for all other fixed points, 
where there is a total of $\chi(M)-(3+2\cdot 1)=19$ (cf. \cite{KS}), the 
local representation is $(z_1,z_2)\mapsto (\mu_p^k z_1,\mu_p^{-k} z_2)$
for some $k\neq 0\bmod p$. To see this, note that the type $\tilde{A}_4$
component $\Lambda_\alpha$ contributes at least $5$ to 
$c_1(K)\cdot [\omega]<7$, so that there is at most one more component
which must be of type (A), i.e., consisting of a single
$J$-holomorphic curve. The claim follows from Proposition 3.7 (1) 
immediately. 

To continue we recall that the signature defect $\text{def}_m$ for an
isolated fixed point $m$ with local representation $(z_1,z_2)\mapsto
(\mu_p^k z_1,\mu_p^{kq} z_2)$ for some $k\neq 0\bmod p$ is 
$$
I_{p,q}\equiv\sum_{k=1}^{p-1}\frac{(1+\mu_p^k)(1+\mu_p^{kq})}
{(1-\mu_p^k)(1-\mu_p^{kq})}.
$$
A direct calculation (cf. Appendix: Proof of Lemma 3.8) gives
$$
I_{5,-1}=4,\;  I_{5,1}=-4,\; I_{5,2}=0,\; I_{5,3}=0.
$$
(Observe the relation $I_{p,q}=-I_{p,-q}$.) We reach a contradiction for 
$p=5$ through the $G$-signature theorem:
$$
-16(5-1)=19 I_{5,-1}+I_{5,1}+I_{5,2} +I_{5,3}+\frac{5^2-1}{3}\cdot (-2)
$$
or equivalently $-64=56$.

It remains to consider for case (ii) the possibility of having a component
$\Lambda_\alpha$ which is represented by a type $\tilde{D}_4$ graph. In
this situation, note that $\Lambda_\alpha$ must be the only 
component (notice that the center vertex in the corresponding graph has
weight $2$, cf. \cite{BPV}, Lemma 2.12 (ii) on page 20, so that 
$\Lambda_\alpha$ contributes at least $6$ to $c_1(K)\cdot [\omega]<7$),
and besides the fixed $(-2)$-sphere, it also contains four isolated fixed
points of rotation numbers $(1,p-2)$. The rest of the fixed points are 
isolated, all of local representation contained in $SL_2(\C)$. By Lemma 2.4,
the $G$-signature theorem implies 
$$
-16 (p-1)\geq 4\cdot I_{p,-2}-\frac{2}{3}(p^2-1),
$$
which is a contradiction because $I_{p,-2}=-I_{p,2}=\frac{1}{6}(p-1)(p-5)$
(cf. Appendix: Proof of Lemma 3.8). 

We have thus shown that there are no nontrivial homologically trivial
finite group actions on $M$ which preserve the symplectic structure 
$\omega$. 

Finally, we remark that there are indeed examples of such $4$-manifolds 
$(M,\omega)$, which actually have nontrivial canonical class $c_1(K)$. 
(This is the case which is not covered in Theorem A.) In fact,
such an $M$ may be obtained by the knot surgery construction of 
Fintushel and Stern \cite{FS1}. More precisely, consider the $K3$ surface
$S_4$ which is the hypersurface in $\P^3$:
$$
S_4\equiv \{[z_0,z_1,z_2,z_3]\in\P^3 \mid z_0^4-z_1^4+z_2^4-z_3^4=0\}
$$
There is a holomorphic elliptic fibration $\pi:S_4\rightarrow \P^1$
such that a generic fiber $F$ of $\pi$ is a cubic curve in a hyperplane 
in $\P^3$ (cf. \S 3.2 of \cite{GS}). Let $\omega_0$ be the K\"{a}hler
form on $S_4$ which is obtained by restricting the Fubini-Study form
to $S_4$. Then $[\omega_0]\in H^2(S_4;\R)$ is an integral class, 
and $[\omega_0]\cdot F=3$ since $F$ is a cubic curve. 

Now applying the knot surgery construction of Fintushel and Stern \cite{FS1}
to $S_4$ at a regular fiber of the holomorphic elliptic fibration with
the K\"{a}hler form $\omega_0$ and using the trefoil knot, one obtains a
symplectic $4$-manifold $(M,\omega)$, such that (1) $M$ is homeomorphic
to $S_4$, (2) the canonical class $c_1(K)$ is the Poincar\'{e} dual of $2F$
where $F$ is the fiber class of the holomorphic elliptic fibration. It
follows easily that $(M,\omega)$ satisfies the required conditions. 

We wish to point out that the $4$-manifold $M$ constructed above is not 
a complex surface (with either 
orientation). To see this, first note that $M$ is not diffeomorphic to 
any $4$-manifold obtained by performing a log transform to a $K3$ surface 
(cf. \cite{FS0, FS1}). In particular, $M$ can not be a complex surface 
with Kodaira dimension $\leq 1$ (cf. e.g. \cite{GS}, Lemma 3.3.4 and
Theorem 3.4.12). On the other hand, $M$ can not be a complex surface of 
general type (necessarily with the opposite orientation) because it violates 
the Miyaoka-Yau inequality $c_1^2\leq 3 c_2$. Hence $M$ is not
a complex surface. Note that this in particular shows that the 
homological rigidity obtained in this example goes beyond that of holomorphic 
actions on K\"{a}hler surfaces in Peters \cite{Pe2}.

\end{example}

%\newpage

\centerline{\bf Appendix: Proof of Lemma 3.8}

\vspace{4mm}

Recall that (cf. Theorem 2.3) for an isolated fixed point $m\in M^G$, 
the signature defect 
$\text{def}_m$ is given by the following expression if the 
local representation at $m$ is 
$(z_1,z_2)\mapsto (\mu_p^k z_1,\mu_p^{kq} z_2)$ for some $k\neq 0\bmod p$
and $q\neq 0\bmod p$:
$$
I_{p,q}\equiv \sum_{k=1}^{p-1}\frac{(1+\mu_p^k)(1+\mu_p^{kq})}
{(1-\mu_p^k)(1-\mu_p^{kq})}
$$
With this notation, it is easily seen that
$$
\text{def}_{(1)}=I_{p,-1}, \text{def}_{(2)}=I_{p,-6}+I_{p,(p+3)/2},
\text{def}_{(3)}=I_{p,2}+2 I_{p,-4}, \text{def}_{(4)}=
I_{p,1}+3 I_{p,-3}.
$$
By Lemma 2.4, $\text{def}_{(1)}=\frac{1}{3}(p-1)(p-2)$ and $I_{p,1}=-I_{p,-1}
=-\frac{1}{3}(p-1)(p-2)$.

To calculate $\text{def}_{(k)}$ for $k=2,3$ and $4$, we go back to the proof
of Lemma 2.4. Recall that $I_{p,q}=-4p\cdot s(q,p)$, and the Dedekind sum
$s(q,p)$ can be computed from
$$
6p\cdot s(q,p)=(p-1)(2pq-q-\frac{3p}{2})-6f_p(q)
$$
where $f_p(q)=\sum_{k=1}^{p-1} k[\frac{kq}{p}]$. Note that 
$I_{p,2}=-\frac{1}{6}(p-1)(p-5)$ follows directly from
$$
f_p(2)=\sum_{k=1}^{p-1} k[\frac{2k}{p}]=\sum_{k=(p+1)/2}^{p-1} k=
\frac{1}{8}(3p-1)(p-1).
$$

To calculate $I_{p,q}$ for the other values of $q$, we recall the following
equation from \cite{HZ} (equation (8) on page 94)
$$
\sum_{k=1}^{p-1}[\frac{kq}{p}]^2-\frac{2q}{p}\cdot \sum_{k=1}^{p-1} 
k[\frac{kq}{p}]=\frac{1}{6p} (1-q^2)(p-1)(2p-1)
$$
and eliminate $f_p(q)$ from the expression for $6p\cdot s(q,p)$.
Consequently, we have
$$
I_{p,q}=(-\frac{2pq}{3}+\frac{q}{3}+\frac{1}{3q}+p-\frac{2p}{3q})(p-1)+
\frac{2p}{q}\cdot \sum_{k=1}^{p-1}[\frac{kq}{p}]^2.
$$

Next we evaluate $\sum_{k=1}^{p-1}[\frac{kq}{p}]^2$ for $q=-4, -6,
\frac{p+3}{2}$ and $q=-3$.

(1) $q=-4$. We shall consider separately the cases $p=4r+1$ and $p=4r+3$.
For $p=4r+1$, we have
$$
\sum_{k=1}^{p-1}[\frac{-4k}{p}]^2=\sum_{k=1}^r (-1)^2+
\sum_{k=r+1}^{2r} (-2)^2+\sum_{k=2r+1}^{3r} (-3)^2+
\sum_{k=3r+1}^{4r} (-4)^2=30r.
$$
For $p=4r+3$, we have
$$
\sum_{k=1}^{p-1}[\frac{-4k}{p}]^2=\sum_{k=1}^r (-1)^2+
\sum_{k=r+1}^{2r+1} (-2)^2+\sum_{k=2r+2}^{3r+2} (-3)^2+
\sum_{k=3r+3}^{4r+2} (-4)^2=30r+13.
$$

(2) $q=-6$. We shall consider separately the cases $p=6r+1$ and $p=6r+5$.
For $p=6r+1$, we have
$$
\sum_{k=1}^{p-1}[\frac{-6k}{p}]^2=\sum_{j=1}^6\sum_{k=(j-1)r}^{jr} (-j)^2
=91 r.
$$
For $p=6r+5$, we have 
$$
\sum_{k=1}^{p-1}[\frac{-6k}{p}]^2=\sum_{k=1}^r (-1)^2+
\sum_{j=1}^4\sum_{k=jr+j}^{(j+1)r+j} (-j-1)^2 +\sum_{k=5r+5}^{6r+4} (-6)^2
=91 r+54.
$$

(3) $q=\frac{p+3}{2}$. Observe that when $k=(2l-1)$ is odd, we have
$$
[\frac{kq}{p}]=\left \{\begin{array}{ll}
l-1 & \mbox{ if } l=1,\cdots, [\frac{p+3}{6}]\\
l & \mbox{ if } l=[\frac{p+3}{6}]+1,\cdots, \frac{p-1}{2},\\
\end{array} \right .
$$
and when $k=2l$ is even, we have 
$$
[\frac{kq}{p}]=\left \{\begin{array}{ll}
l & \mbox{ if } l=1,\cdots, [\frac{p}{3}]\\
l +1 & \mbox{ if } l=[\frac{p}{3}]+1,\cdots, \frac{p-1}{2}.\\
\end{array} \right .
$$
A direct calculation gives 
\begin{eqnarray*}
\sum_{k=1}^{p-1}[\frac{kq}{p}]^2 & = & \sum_{l=1}^{[(p+3)/6]} (l-1)^2+
\sum_{l=[(p+3)/6]+1}^{(p-1)/2} l^2+ \sum_{l=1}^{[p/3]} l^2 +
\sum_{l=[p/3]+1}^{(p-1)/2} (l+1)^2\\
    & = & \left\{\begin{array}{ll}
r(18r^2+13 r+3) & \mbox{ if } p=6r+1\\
(r+1)(18r^2+31 r+14) & \mbox{ if } p=6r+5.\\
\end{array}\right .
\end{eqnarray*}

(4) $q=-3$. We shall consider separately the cases $p=3r+1$ and $p=3r+2$.
For $p=3r+1$, we have
$$
\sum_{k=1}^{p-1}[\frac{-3k}{p}]^2=\sum_{k=1}^r (-1)^2+
\sum_{k=r+1}^{2r} (-2)^2+\sum_{k=2r+1}^{3r} (-3)^2=14 r.
$$
For $p=3r+2$, we have
$$
\sum_{k=1}^{p-1}[\frac{-3k}{p}]^2=\sum_{k=1}^r (-1)^2+
\sum_{k=r+1}^{2r+1} (-2)^2+\sum_{k=2r+2}^{3r+1} (-3)^2=14 r+4.
$$

Now plug these formulas back into the expression for $I_{p,q}$. We obtain
\begin{eqnarray*}
I_{p,-4} & = & \left\{\begin{array}{ll}
\frac{4}{3} r(r-4) & \mbox{ if } p=4r+1\\
\frac{2}{3}(2r^2+1) & \mbox{ if } p=4r+3,\\
\end{array}\right .\\
I_{p,-6} & = & \left\{\begin{array}{ll}
2r(r-6) & \mbox{ if } p=6r+1\\
2r^2+4r+4 & \mbox{ if } p=6r+5,\\
\end{array}\right .\\
I_{p,(p+3)/2} & = & \left\{\begin{array}{ll}
2r(2-r) & \mbox{ if } p=6r+1\\
-2r^2+4r+4 & \mbox{ if } p=6r+5,\\
\end{array}\right .\\
I_{p,-3} & = & \left\{\begin{array}{ll}
r(r-3) & \mbox{ if } p=3r+1\\
r(r-1) & \mbox{ if } p=3r+2.\\
\end{array}\right .
\end{eqnarray*}
This gives the formulas of $\text{def}_{(k)}$ for $k=2,3$ and $4$, and
Lemma 3.8 follows.

%\newpage

{\Small 
W. Chen: Department of Math. and Stat., University of Massachusetts,
Amherst, MA 01003.\\
{\it e-mail:} wchen@math.umass.edu

S. Kwasik: Mathematics Department, Tulane University,
New Orleans, LA 70118. \\
{\it e-mail:} kwasik@math.tulane.edu\\
}

\end{document}